\documentclass[11pt]{article} % use larger type; default would be 10pt

\usepackage[utf8]{inputenc} % set input encoding (not needed with XeLaTeX)

\usepackage{geometry} % to change the page dimensions
\usepackage{graphicx} % support the \includegraphics command and options

\usepackage{booktabs} % for much better looking tables
\usepackage{array} % for better arrays (eg matrices) in maths
\usepackage{paralist} % very flexible & customisable lists (eg. enumerate/itemize, etc.)
\usepackage{verbatim} % adds environment for commenting out blocks of text & for better verbatim
\usepackage{subfig} % make it possible to include more than one captioned figure/table in a single float
\usepackage{amsmath,amsthm,amssymb,color,epic,eepic,cite}
\usepackage{hyperref}
\usepackage{tikz}
\usetikzlibrary{matrix}
\usepackage{fancyhdr} % This should be set AFTER setting up the page geometry
\usepackage{sectsty}
\allsectionsfont{\sffamily\mdseries\upshape} % (See the fntguide.pdf for font help)
\usepackage[nottoc,notlof,notlot]{tocbibind} % Put the bibliography in the ToC
\usepackage[titles,subfigure]{tocloft} % Alter the style of the Table of Contents

\newcommand{\Z}{{\Bbb Z}} %??
\newcommand{\C}{{\Bbb C}} %??
\newcommand{\K}{{\Bbb K}}
\newcommand{\A}{{\cal A}}

\newcommand{\hbs}{{\bar{H}^*}}
\newcommand{\tpt}{{\tilde{\otimes}}} 
\newcommand{\D}{{\mathcal{D}}}

\newcommand{\cD}{{\cal D}}

\newcommand{\B}{{\cal B}}

\renewcommand{\H}{{\cal H}}
\newcommand{\la}{\lambda}

\newcommand{\al}{\alpha}
\newcommand{\e}{\epsilon}
\newcommand{\vep}{\varepsilon}

\newcommand{\hf}{\widehat{f}}

\newcommand{\noi}{{\noindent}}
\newcommand{\nn}{{\nonumber}}
\newcommand{\bea}{\begin{eqnarray}}
\newcommand{\ena}{\end{eqnarray}}
\newcommand{\beit}{\begin{itemize}}
\newcommand{\enit}{\end{itemize}}

\newcommand{\be}{\begin{eqnarray*}}
\newcommand{\en}{\end{eqnarray*}}
\newcommand{\lb}[1]{\label{#1}}
\newcommand{\id}{{\rm id}}

\def\infq4p#1{{(#1;q^4,p)_\infty}}

\newcommand{\tot}{\widetilde{\otimes}}

\newcommand{\mmatrix}[1]{\begin{matrix} #1 \end{matrix}}

\font\teneufm=eufm10
\font\seveneufm=eufm7
\font\fiveeufm=eufm5
\newfam\eufmfam
\textfont\eufmfam=\teneufm
\scriptfont\eufmfam=\seveneufm
\scriptscriptfont\eufmfam=\fiveeufm
\def\frak#1{{\fam\eufmfam\relax#1}}
\let\goth\frak
\newcommand{\slth}{\widehat{\goth{sl}}_2}
\newcommand{\glth}{\widehat{\goth{gl}}_2}
\newcommand{\slt}{\goth{sl}_2}
\newcommand{\glnh}{\widehat{\goth{gl}}_N}

\newcommand{\h}{H}
\font\seventeeneufm=eufm10 scaled\magstep3   %for section
   
\newcommand{\slthBig}{\widehat{\mbox{\seventeeneufm sl}}_2} %??
\newcommand{\glthBig}{\widehat{\mbox{\seventeeneufm gl}}_2} %??
\newcommand{\fp}{{F^+(u_1)}}

\newcommand{\fpm}{{F^\pm(u_1)}}
\newcommand{\fmp}{{F^\mp(u_1)}}
\newcommand{\Ep}{{E^+(u_1)}}

\newcommand{\fpt}{{F^+(u_2)}}

\newcommand{\fpmt}{{F^\pm(u_2)}}
\newcommand{\fmpt}{{F^\mp(u_2)}}
\newcommand{\Ept}{{E^+(u_2)}}
\newcommand{\Emt}{{E^-(u_2)}}

\newcommand{\kpone}{{K_1^+(u_1)}}

\newcommand{\kpmtwo}{{K_2^\pm(u_1)}}
\newcommand{\kmptwo}{{K_2^\mp(u_1)}}
\newcommand{\kponet}{{K_1^+(u_2)}}

\newcommand{\kmtwot}{{K_2^-(u_2)}}
\newcommand{\kpmtwot}{{K_2^\pm(u_2)}}

\makeatletter
\@addtoreset{equation}{section}
\makeatother

\newtheorem{thm}{Theorem}[section]
\newtheorem{prop}[thm]{Proposition}
\newtheorem{lem}[thm]{Lemma}
\newtheorem{cor}[thm]{Corollary}

\newtheorem{dfn}[thm]{Definition}

\title{Dynamical FRT construction of $U_{q,x}(gl_N)$}
\author{Bharath Narayanan\\
The Pennsylvania State University\\
narayana@math.psu.edu}
\begin{document}
\bibliographystyle{unsrt}

%\vspace{2cm}
\begin{center}
%{\Large \bf Two Equivalent realizations of the
% Dynamical Affine Quantum Group
% $U_{q,x}(\glthBig)=F_{q,\la}[\widehat{GL_2}]$ .\\[10mm] }
{\Large \bf Two Equivalent Realizations of Trigonometric Dynamical Affine Quantum Group $U_{q,x}(\widehat{sl_2})=U_{q,\lambda}(\widehat{sl_2})$, Drinfeld Currents and Hopf Algebroid Structures\\[10mm]}
{\large  Bharath Narayanan}\\[6mm]
{\it Department of Mathematics, 
\\The Pennsylvania State University\\
       narayana@math.psu.edu}\\[10mm]
\end{center}

\begin{abstract}
\noindent 
Two new realizations, denoted $U_{q,x}(\widehat{gl_2})$ and $U(R_{q,x}(\widehat{gl_2}))$ of the dynamical quantum affine algebra $U_{q,\lambda}(\widehat{gl_2})$ are proposed, based on Drinfeld-currents and $RLL$ relations respectively, along with a Heisenberg algebra $\left\{P,Q\right\}$, with $x=q^{2P}$. Here $P$ plays the role of the dynamical variable $\lambda$ and $Q=\frac{\partial}{\partial P}$. An explicit isomorphism from $U_{q,x}(\widehat{gl_2})$ to $U(R_{q,x}(\widehat{gl_2}))$ is established, which is a dynamical extension of the Ding-Frenkel isomorphism of $U_{q}(\widehat{gl_2})$ with $U(R_{q}(\widehat{gl_2}))$ between the Drinfeld realization and the Reshetikhin-Tian-Shanksy construction of quantum affine algebras. Hopf algebroid structures and an affine dynamical determinant element are introduced and it is shown that $U_{q,x}(\widehat{sl_2})$ is isomorphic to $U(R_{q,x}(\widehat{sl_2}))$. The dynamical construction is based on the degeneration of the elliptic quantum algebra $U_{q,p}(\widehat{sl_2})$ of Jimbo, Konno et al. as the elliptic variable $p \to 0$.
\end{abstract}

\section{Introduction}

The {\em elliptic} affine quantum group $U_{q,p}(\slth)$, where $q$ denotes the quantum variable and p the elliptic, is studied in detail by M. Jimbo,  H. Konno et al. in \cite{Konno, JKOS2} using elliptic deformations (i.e. twists) of the quantum Drinfeld currents and a Heisenberg algebra $\H$ containing the dynamical variable $P$ and a dual element $Q$,  to construct an operator $L^+(u)$ which obeys the $RLL$ relations:
\bea \lb{ellipticrll}
&&R^{+(12)}(u_1-u_2,P+h)L^{+ (1)}\left(u_1\right)L^{+ (2)}\left(u_2\right)=L^{+ (2)}\left(u_2\right)L^{+ (1)}\left(u_1\right)R^{+*(12)}(u_1-u_2,P)\nn\\
\ena

\noi
 associated to the {\em {elliptic R-matrix}}:
\bea
&&R_{q,p}^+(u,P)=\rho^+(u)
\left(
\begin{array}{cccc}
1 &0                  &0             &0 \\
0  &\frac{[P+1] [P-1][u]  }{[P]^2[1+u]}    &\frac{[1][P+u]}{[P][1+u]} &0  \\
 0 &\frac{[1][P-u]}{[P][1+u]}&\frac{[u]}{[1+u]}&0 \\
0  & 0                 &0             &1 \\
\end{array}
\right) ,
\ena
\lb{Rmat4}

\noi
where  $\rho^+(u)$ is a suitably chosen coefficient, $R^{+*}(u,P)=R^{+}(u,P)|_{r\to r^*}$ with $r^*=r-c$
and the Jacobi theta function is given by
\be
&&[u]=\frac{q^{\frac{u^2}{r}-u}}{(p;p)_\infty^3}\Theta_p(z),\qquad p=q^{2r}, \qquad z=q^{2u},\\
&&\Theta_p(z)=(z;p)_\infty(p/z;p)_\infty(p;p)_\infty, \qquad (z;p)_\infty=\prod_{n=0}^\infty(1-zp^{n}).
\en
 
\begin{figure}[here]
\includegraphics[width=1.0\textwidth]{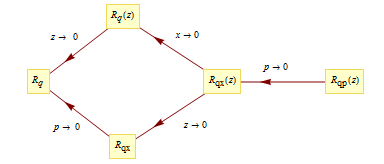}
\caption{$R$-matrix degenerations}
\label{fig:figure1}
\end{figure}

Most of the degenerations of this $R$-matrix, denoted $R_{q,p}(z)$, yield well-known $R$-matrices (see Figure \ref{fig:figure1}).

The top entry is the quantum affine $R$-matrix $R_{q}(z)=R_q(\glth)$ \eqref{quantumrmatrix} considered by Reshetikhin and Semenov-Tian-Shansky in \cite{RS}.  The lowest entry is the quantum dynamical $R$-matrix  $R_{q,x}(\slt)$ obtained by a twist construction using the quantum $6j$-symbols \cite{BBB, ESS, KR}.
  The left-most $R$-matrix $R_q$ corresponds to the standard Drinfeld-Jimbo quantum group.  The second one from the right, $R_{qx}(z)=R_{q,x}(u,P)$ (see \eqref{rplus}) is the $R$-matrix considered in this paper. Here $x=q^{2P}$, where $P$ is a generator of the Heisenberg algebra $\H=\left\{Q,P\right\}$ and plays the role of the dynamical variable $\la$ in the more standard formulations \cite{EV1, EV2, KR, EF}, while $Q=\frac{\partial}{\partial P}$.  The exact relation of $R_{qx}$ to the $R$-matrix $R(\lambda)$ in \cite{KR} is given by $R(\lambda)=\lim_{z \to 0}R^{(21)}(u,P) \mid_{ q= q^{-1}}$.

Based on each of these $R$-matrices,  one can define a pair of algebras $(U(R), A(R))$, each of which carries a Hopf-type structure and can be considered as appropriately defined twists (or deformations) of the underlying universal enveloping algebra $U(\slt)$ and coordinate function algebra $\C[SL_2]$, respectively.

\begin{figure}[here]
\includegraphics[width=1.05\textwidth]{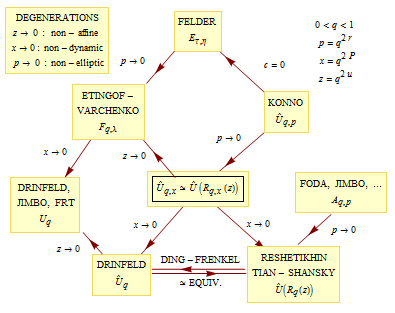}%newestuq.png}
\caption{Dynamical Algebra Degenerations for $sl_2$}
\label{fig:algebradegenerations}
\end{figure}

  The algebra of primary concern in this article will be $U(R)$ for $R=R_{q,x}(u,P)$, located in the center of Figure \ref{fig:algebradegenerations}.

Jimbo, Konno et al.\cite{JKOS2} use elliptic Drinfeld currents to define the algebra  $\widehat{U}_{q,p}=U_{q,p}(\slth)$ which can be viewed as a tensor product of the underlying quantum affine algebra, $\widehat{U}_{q}=U_{q}(\slth)$ and a Heisenberg algebra which includes the elliptic variable $p$.  They use only {\em positive} half-currents to define an $L$-operator that satsfies the $RLL$ relation \eqref{ellipticrll}.

There is another construction of the elliptic affine algebra, denoted $E_{\tau,\eta}(\slth)$, by Enriquez and Felder \cite{EF}, using both positive and negative half-currents (for the precise relation to $U_{q,p}(\slth)$ see Section (6.2) in \cite{JKOS2}).
At $c=0$, we can identify $U_{q,p}(\slth)$ with Felder's original elliptic quantum group  $E_{\tau,\eta}=E_{\tau,\eta}[SL_2]$ (the topmost box) \cite{Felder}.  A vertex-type non-dynanical elliptic quantum algebra $A_{q,p}=A_{q,p}(\slth)$ is investigated in \cite{yan}.

J. Ding and I. Frenkel \cite{DF} proved an equivalence $U_{q}(\glnh) \simeq U(R_q(\glnh))$, between two realizations of the quantum affine algebra associated with $\glnh$,  in terms of Drinfeld currents \cite{Drinfeld} and via $RLL$-relations \cite{RS}, due to Reshetikhin and Semenov Tian-Shanksy.   The non-elliptic limit of $A_{q,p}(\slth)$ is also considered in \cite{yan} and it is speculated there that it must be the same as the quantum affine algebra $U(R_q(\slth))$.  The dynamical quantum group $F_{q,\la}=F_{q,\la}[SL_2]$ is discussed in \cite {EV1, KR}, based on an FRST-construction.   Farthest on the left in Figure \ref{fig:algebradegenerations} sits $U_q=U_q(\slt)$, the usual Drinfeld-Jimbo quantum group.

Clearly, only the box that is in the center of Figure \ref{fig:algebradegenerations} remains to be investigated.  Since it is only one link away from four other boxes, it is natural to expect that it has all the desirable characteristics of those algebras.   The purpose of this article is to confirm this expectation by defining two dynamical quantum algebras  $U_{q,x}(\glth)$ and $U(R_{q,x}(\glth))$,  endowing them with suitable H-bialgebroid structures and proving that they are isomorphic as $H$-algebras.  Thus we extend the corresponding result of Ding-Frenkel to the dynamical case (for $N=2$).  We also confirm that the generators and relations for the degenerations in the non-dynamical and non-affine directions are consistent with the known models.  As mentioned earlier, one recovers $U_{q,x}(\slth)$ in the limit as $p \to 0$ of the elliptic algebra $U_{q,p}(\slth)$  \cite {JKOS,K2}, allowing us to identify $U_{q,x}(\slth)$ with $U_{q,p}(\slth)/\ p \cdot U_{q,p}(\slth)$.  Our algebra is similar, but not identical, to the construction by P. Xu in \cite{Ping} of  the dynamical quantum groupoid for {\em finite-dimensional} Lie algebras.

It turns out that the key to transit most efficiently from the quantum to dynamical quantum world is the introduction, by H. Konno \cite{Konno} of the Heisenberg algebra $H=\left\{P, Q\right\}$.  %The elliptic algebra can be viewed as the tensor product of the usual quantum affine algebra $U_q(\slth)$ with $\H$. 
Then $U_{q,x}(\slth)$ can be viewed as a semidirect product (smash product) algebra   
$\left(U_{q}(\slth) \otimes M_{H^*} \right) \otimes \C[H_Q^*]$.
%$U_{q}(\glth) \otimes \C[H]$ \# $\C[H_Q^*]$.  Here $H=\C P$ and $H_Q=\Z Q$. 
  An important difference from the elliptic case is the observation that, at $p=0$, the positive half currents $\left\{E^+ (u), F^+(u), K_1^+(u), K_2^+(u)\right\}$ contain only {\em non-positive} powers of $z$.  Therefore,  to recover the total algebra, i.e. non-negative powers, one must include the negative half-currents $\left\{E^- (u), F^-(u), K_1^-(u), K_2^-(u)\right\}$ as well and the full set of %$R_\pm^\pm L _\pm^\pm L _\pm^\pm$ relations,
$RLL$ relations:
\bea \lb{drll}
&&R^{\pm(12)}(u,P+h)L^{\pm (1)}\left(u_1\right)L^{\pm (2)}\left(u_2\right)=L^{\pm (2)}\left(u_2\right)L^{\pm (1)}\left(u_1\right)R^{\pm(12)}(u,P), \nn\\\
&&R^{\pm(12)}\left(u\pm\frac{c}{2},P+h\right)L^{\pm(1)}\left(u_1\right)L^{\mp(2)}\left(u_2\right)=L^{\mp(2)}\left(u_2\right)L^{\pm(1)}\left(u_1\right)R^{\pm(12)}\left(u\mp\frac{c}{2},P\right),\nn\\
\ena
where $R^-(u,P)$ and $R^+(u,P)$ differ by the coefficient $\rho^\pm(z)$ given in \eqref{rho}.  The operators $L^{\pm}(u)$ can be viewed as part of a single $L$-operator which is a  doubly infinite series \cite{yan} with a single RLL relation, but we will not pursue this approach here.
% but this formulation is not well-suited for calculations in terms of the Drinfeld half-currents and for its representation theory. 
We will also consider the algebras $U_{q,x}^{+}(\slth)$ (resp. $U_{q,x}^{-}(\slth)$)  defined by the equations containing $L^{+}$ (resp $L^{-}$) only, consisting of non-positive (resp. non-negative) powers of $z$ only. 

The main results proven in this article are:
\begin{enumerate}
\item $U_{q,x}(\glth) \simeq U(R_{q,x}(\glth))$ as $H$-algebras.
\item The subalgebras $U_{q,x}^{\pm}(\glth)$ and $U(R^\pm_{q,x}(\glth))$ are $H$-Hopf algebroids.  
\item The total algebras $U_{q,x}(\glth)$ and $U(R_{q,x}(\glth)$ are $H$-Hopf algebroids.
\item The results (i) - (iii) remain true upon replacing $\glth$ by $\slth$ .
\end{enumerate}

\noi
 (see Theorems \ref {mainthm}, \ref{thmUR}, \ref{mainthmsl2} and  \ref{thmuq} ).

%It turns out that there is a natural inherent symmetry between the plus and minus.
Let us mention some motivation for the current investigation.   In \cite{clad}, Arnaudon et al. illustrate the various degenerations and their inverses (twists) of the elliptic algebras and double Yangians.  They derive the $R$-matrix and write the basic $RLL$ relations for $U_{q,\la}(\slth)$ but do not perform the Lax expansion.   They also stress the physical importance, relevant to the Calegoro-Moser systems, and point out that the "dynamical" variable $\la$ (which is elevated by Jimbo, Konno et al. in \cite{JKOS,Konno} to an actual generator $P$ for the algebra $U_{q,p}(\slth)$) can be identified with the momentum of the system.  The results in the current article provide a positive step towards confirming  their expectation that similar genuinely dynamical structures exist for all formal limits (or twists) of the quantum algebras described there, and play an important role in solving the models where such algebras arise (see the Conclusion in \cite{clad}).  This article creates the framework to develop the representation theory of $U_{q,x}(\slth)$. 
% The first motivation for was to describe the irreducible modules of twisted quantum algebras, based on a question asked to the author by P. Etingof.
The companion article \cite{merep} examines the finite dimensional representations of $U_{q,x}(\slth)$ and its relation to representations of $U_{q,p}(\slth)$ and hypergeometric series.
The main algebra introduced in this paper, and its representations, can be considered part of the {\em algebraic analysis} framework of solving statistical models by using the representation theory of infinite dimensional algebras, due to M. Jimbo.   Specifically,  we are in the Andrews-Baxter-Forrester regime, studying (degenerations of) solutions of the DQYBE for the 8-vertex and RSOS models.  A deep result that the DQYBE is equivalent to the Star-Triangle relation was shown by Felder in \cite{Felder}.  The relation of $U_{q,x}(\slth)$ with the Universal vertex-IRF transformation  is given by Buffenoir et al. in \cite{BRT}.  Further, the scaling limit as $q \to 1$ of the elliptic quantum algebra is  studied by Hou and Yang \cite{Hou}.
\vspace{-2pt}

{\em Outline of the Article}.   Section 2 reviews the Ding-Frenkel construction of the quantum affine algebra isomorphism $U_q(\glth)\simeq U(R_q(\glth))$.   In Section 3, the $H$-algebra $U_{q,x}(\glth)$ is defined using Drinfeld currents and a Heisenberg algebra.  Further, the dynamical half-currents are defined and their commutation relations are proven.  The $H$-algebra $U(R_{q,x}(\glth))$ is defined in Section 4 where we prove our main result, that $U_{q,x}(\glth) \simeq U(R_{q,x}(\glth))$.  Section 5 is devoted to $H$-Hopf-algebroid structures,  the dynamical determinant element  and the subalgebra $U_{q,x}(\slth)$.  The basic definitions are included at the beginning of the section.  The expansions of the RLL relations appear in Appendix A. 
%\newpage
\section {The Quantum Affine Algebra \texorpdfstring{$U_{q}(\glthBig)$}{affine gl2}}
In this section we review the definition of the Drinfeld realization of the quantum affine algebra $U_{q}(\slth)$, as given in \cite{DF}.   The quantum affine algebra $U_{q}(\slth)$ has been defined in terms of Chevalley generators (Kac-Moody algebras), Drinfeld currents (Yangians and quantum loop algebras), and Reshetekhin-Tian-Shansky's  $RLL$-algebra (FRT-type construction).  We consider the last two models here.

\subsection{Definition of \texorpdfstring{$U_{q}(\glthBig)$}{affine Uq}}\lb{app:1.3}
Let us recall the results of \cite{DF} which are relevant to this treatise.  The definition of the Drinfeld realization of the quantum affine algebra $U_q(\glth)$ is adapted from Ding-Frenkel's definition \cite{DF,FJ} of  $U_q(\glnh)$.{\footnote{The relation with the more standard presentation of quantum affine $\slth$ is given in Section \ref{drinfeld}}}.  Let $q$ be a complex number $q\not=0$ such that $|q|<1$. 
\begin{dfn}\cite{Drinfeld}\lb{defUq}
For a field $\K \supseteq \C$, the quantum affine algebra $\K[U_q(\glth)]$ 
in the Drinfeld realization is an associative algebra over $\K$ 
generated by the generators $q^{\pm h}$,  $q^{\pm c}$, $k_1^{\pm}(z), k_2^{\pm}(z)$, $e_n, f_n (n\in \Z)$.    
The defining relations are given as follows:
\be
&&q^c :\hbox{ central },\nn\\
&&q^hq^{-h}=q^{-h}q^h=1,\\
&& k_1^{+ }(z)k_1^{-}(w)=k_1^{- }(w)k_1^{+ }(z), \quad k_2^{+}(z)k_2^{- }(w)=k_2^{- }(w)k_2^{+ }(z), \nn\\
&& k_1^{\pm }(z)k_1^{\pm }(w)=k_1^{\pm }(w)k_1^{\pm }(z), \quad k_2^{\pm }(z)k_2^{\pm }(w)=k_1^{\pm }(w)k_1^{\pm }(z), \nn\\
%&&\frac{{zq}^{\frac{c}{2}}-{wq}^{-\frac{c}{2}}}{{zq}^{-\frac{c}{2}-1}-  {wq}^{-\frac{c}{2}}q}k_1^-(z)k_2^+(w)=k_2^+(w)k_1^-(z)\frac{  {zq}^{-\frac{c}{2}}-  {wq}^{\frac{c}{2}}}{  {zq}^{-\frac{c}{2}-1}-  {wq}^{\frac{c}{2}+1}}, \nn\\
&&\frac{ {zq}^{-\frac{c}{2}}- {wq}^{\frac{c}{2}}}{ {zq}^{-\frac{c}{2}-1}- {wq}^{\frac{c}{2}+1}}k_1^+(z)k_2^-(w)=k_2^-(w)k_1^+(z)\frac{ {zq}^{\frac{c}{2}}- {wq}^{-\frac{c}{2}}}{ {zq}^{\frac{c}{2}-1}- {wq}^{-\frac{c}{2}+1}} \hspace{0.1cm}, \nn\\
&&\frac{{zq}^{\frac{c}{2}}-{wq}^{-\frac{c}{2}}}{{zq}^{\frac{c}{2}+1}-  {wq}^{-\frac{c}{2}-1}}k_2^+(z)k_1^-(w)=k_1^-(w)k_2^+(z)\frac{  {zq}^{-\frac{c}{2}}-  {wq}^{\frac{c}{2}}}{  {zq}^{-\frac{c}{2}+1}-  {wq}^{\frac{c}{2}-1}}\hspace{0.1cm} , \nn\\
&&k_1^{\pm }(z)e(w)k_1^{\pm }(z)^{-1}=\frac{ {zq}^{\mp \frac{c}{2}}-w}{ {zq}^{\mp \frac{c}{2}-1}- {wq}}e(w), \quad
k_2^{\pm }(z)e(w)k_2^{\pm }(z)^{-1}=\frac{ {zq}^{\mp \frac{c}{2}}-w}{ {zq}^{\mp \frac{c}{2}+1}- {wq}^{-1}}e(w), \nn\\
&&k_1^{\pm }(z)^{-1}f(w)k_1^{\pm }(z)=\frac{ {zq}^{\pm \frac{c}{2}}-w}{ {zq}^{\pm \frac{c}{2}-1}- {wq}}f(w), \quad
k_2^{\pm }(z)^{-1}f(w)k_1^{\pm }(z)=\frac{ {zq}^{\pm \frac{c}{2}}-w}{ {zq}^{\pm \frac{c}{2}+1}- {wq}^{-1}}f(w), \nn\\
&&(zq^{-1}-wq)
e(z)e(w)= (zq-wq^{-1}) e(w)e(z), \quad
(zq-wq^{-1})
f(z)f(w)= (zq^{-1}-wq) f(w)f(z),\nn
\\
%next eqution different from theirs by c to -c
&&[e(z),f(w)]=(q-q^{-1})
\left(\delta\left(q^{c}\frac{z}{w}\right)\psi(q^{\frac{-c}{2}}w)
-\delta\left(q^{-c}\frac{z}{w}\right)\varphi(q^{\frac{c}{2}}w)
\right)\nn
\lb{uqslth},\\
%en
%with  
%\be
&&\hspace{-5pt} {\text { where }}\ e(z)=\sum_{n\in \Z}e_{n} z^{-n},\quad
f(z)=\sum_{n\in \Z}f_{n} z^{-n},\quad
[n]_q=\frac{q^n-q^{-n}}{q-q^{-1}},\\
&&\quad \qquad \psi(z)=\sum_{n \ge 0}\psi_{n} z^{-n},\quad
\varphi(z)=\sum_{n \ge 0}\varphi_{-n} z^{n},\quad
\delta(z)=\sum_{n\in\Z}z^n,\\
&&\quad \qquad \psi(z)=k_1^{+}(z)k_2^{+}(z)^{-1},\quad 
\varphi(z)=
k_1^{-}(z)k_2^{-}(z)^{-1}.
\en
\end{dfn}
Note that $\psi (0)=q^h$,  $\varphi (0)=q^{-h}$ and the final defining relation in Definition \ref {defUq} can be equivalently stated in terms of the generating modes as:
\bea \lb{efcomm}
&&[e_m,f_n] =(q-q^{-1})\left(q^{\frac{c(n-m)}{2}}\psi _{m+n}-q^{\frac{c(m-n)}{2}}\varphi _{m+n}\right) \qquad m, n \in \Z.
\ena

\subsection{Ding-Frenkel's Equivalence} \lb{subsecdf}
We summarize the results in \cite{DF} on the quantum affine algebra  $U_q(\widehat{gl_n})$ at $n=2$.
For $0<q<1$, consider the quantum $R$-matrix:
\begin{center}
$R_q(z)= \lb{quantumrmatrix}
\left(
\begin{array}{cccc}
 1 & 0 & 0 & 0 \\
 0 & \frac{q(1-z)}{1-q^2z} & \frac{1-q^2}{1-q^2z} & 0 \\
 0 & \frac{z\left(1-q^2\right)}{1-q^2z} & \frac{q(1-z)}{1-q^2z} & 0 \\
 0 & 0 & 0 & 1 \\
\end{array}
\right).$
\end{center}
\noi
Reshetikhin and Semenov-Tian-Shanksy \cite{RS} define the quantum affine algebra $U(R_q)$ by: 
\begin{dfn} \lb{RSdefn}
 $U(R_q)$ is an associative algebra over $\C$ with central element $c$,  generated by $L^{\pm}(z)=\left(L^\pm_{ab}(z)\right)_{a,b=1}^2$, with $L^\pm_{ab}(z)=\sum_{n=0}^\infty L_{ab,\pm n}z^{\mp n}$ such that $L_{aa,o}^+ L_{aa,o}^- = L_{aa,o}^- L_{aa,o}^+ =1$ and the affine RLL-relations:
\be
&&R_q\left(\frac{z}{w}\right)L_1^{\pm }(z)L_2^{\pm }(w)=L_2^{\pm }(w)L_1^{\pm }(z)R_q\left(\frac{z}{w}\right), \\
&&R_q\left(q^c\frac{z}{w}\right)L_1^+(z)L_2^-(w)=L_2^-(w)L_1^+(z)R_q\left(q^{-c}\frac{z}{w}\right),
\en
where we denote
\be
L_1^\pm(z)=L^\pm(z) \otimes 1 \qquad \text {and} \qquad
L_2^\pm(z)=1 \otimes L^\pm(z).
\en
\end{dfn}
 Some extra relations involving an auxillary operator $\tilde{L}(z)$,  appearing in the original definition \cite{RS}  are absorbed into Definition \ref{RSdefn} by imposing that $\tilde{L}(z)$ satisfies the following natural condition (see eq(3.21) in \cite{DF}):
\be
&& \tilde{L}(z)=((L^{\pm}(z))^t)^{-1}.
\en
\noi
It is known \cite{DF} that $L^\pm(z)$ have the following unique Gauss decompositions:
%\footnote{For the exact correspondence to \cite{DF} use $L^{-1}$, $(k_{i}^{\pm})^{-1}$ and $-c$ in place of $L$, $k_i^{\pm}$ and $c$.}
\be
&&L^{\pm }(z)=\left(
\begin{array}{cc}
 1 & f^{\pm }(z) \\
 0 & 1 \\
\end{array}
\right)\left(
\begin{array}{cc}
 k_1^{\pm }(z) & 0 \\
 0 & k_2^{\pm }(z) \\
\end{array}
\right)\left(
\begin{array}{cc}
 1 & 0 \\
 e^{\pm }(z) & 1 \\
\end{array}
\right),
\en
where  $e^{\pm }(z), f^{\pm }(z)$, $k_1^{\pm }(z)$ and $k_2^{\pm }(z)$ are elements in $U(R_q)$, with $k_i^{\pm }(z)$ invertible.  
\noi
Ding-Frenkel's main result guarantees that the Drinfeld realization is isomorphic to Reshetikhin and Semenov-Tian-Shanksy's presentation:
\vspace{-1pt}
\begin{thm} \lb{dfthm} \cite{DF}
There exists an algebra isomorphism
\vspace{-6pt}
\be
&&\phi : U_{q}(\glth) \longrightarrow U(R_q),\\
&&\hspace{32pt} e(z) \mapsto  e^+\left(q^{-c/2}z\right)-e^-\left(q^{c/2}z\right),\\
&&\hspace{32pt} f(z) \mapsto  f^+\left(q^{c/2}z\right)-f^-\left(q^{-c/2}z\right),\\
&&\hspace{26pt} k_i^{\pm }(z) \mapsto  k_i^{\pm }(z) \text{ and } q^{\pm c }\mapsto  q^{\pm c }.
\en
\end{thm}

\section{The $H$-algebra \texorpdfstring{$U_{q,x}(\glth)$}{Dynamical Affine Uqx(gl2)}}
\subsection{$H$-algebras}  
We will need with the notions of $H$-algebras $(A,H,\mu_l,\mu_r)$, $H$-bialgebroids $(A,H,\mu_l,\mu_r,\Delta,\varepsilon)$ and $H$-Hopf algebroids $(A,H,\mu_l,\mu_r,\Delta,\varepsilon,S)$.  We also need their dynamical tensor products $A \tilde \otimes B$. The definitions are adapted from \cite{bohm,K2,R}. %and their dynamical representations
It is worthwhile to observe that, when $H=0$, an $H$-Hopf algebroid becomes an ordinary Hopf algebra.

\begin{dfn} [$H$-Algebra] \lb{halg}
An associative algebra $A$ over $\C$ is an $H$-algebra if it has an $H$-bigrading such that $A=\oplus_{\alpha,\beta \in H^*}{A_{\alpha,\beta}}$, along with the left and right moment maps, $\mu_l,\mu_r:M_{H^*} \rightarrow A_{00}$ satisfying:
\be
\mu_l(\hf)x=x \mu_l(T_\al \hf), \quad \mu_r(\hf)a=a \mu_r(T_\beta \hf), \qquad 
a\in A_{\al\beta},\ \hf\in M_{\h^*},
\en
where $T_\al$ denotes the automorphism $(T_\al \hf)(\la)=\hf(\la+\al)$ of $M_{\h^*}$.

\end{dfn}
Consider two $H$-algebras $A$ and $B$.  An {\em {H-algebra homomorphism}} $\Psi$ between $A$ and $B$ is an algebra homomorphism which preserves the bigrading and moment maps:
\be
\Psi(A_{\al\beta})\subseteq B_{\al\beta} \  {\text { for all }}\  \alpha,\beta\in \h^* \ {\text {and}}\\
\Psi(\mu^A_l(\hf))=\mu^B_l(\hf), \ \Psi(\mu^A_r(\hf))=\mu^B_r(\hf).
\en
\subsection{The Construction of  \texorpdfstring{$U_{q,x}(\glthBig)$}{Uqx}}
In the dynamical case, there are two main constructions known as the {\em face} and {\em vertex} models which are referred to in the literature as the $\B$ and $\A$ models, respectively.  These algebras are both $U_{q}(\slth)$ as algebras, 
obtained via quasi-Hopf twists (twistors) as described explicitly in \cite{JKOS}.   We will employ the Heisenberg algebra and trignometric Drinfeld currents to construct the dynamical half-currents and define the dynamical quantum affine algebra, $U(R_{q,x}(\glth))$, which is of {\em face} type.% This lays the foundation for defining $U(R_{q,x}(\glth))$ and for proving that $U_{q,x}(\glth) \simeq U(R_{q,x}(\glth))$, which is carried out in the next section. 
\subsubsection{Heisenberg Algebra $\mathcal{H}$ } \lb{heis}
Let us define a Heisenberg algebra $\mathcal{H}$ with generators $P$ and $Q$ such that $[P,Q]=-1$. Denote $H=\mathbb{C}P$, $H^*=\mathbb{C}Q$, $\mathcal{H}=H\oplus H^*$ with the pairing given by $<Q,P> = 1$ and $<x,y> =0$ for all other $x,y$ (for example, we can choose $Q=\frac{\partial }{\partial P}$).
 Let $\hbs = \text{$\mathbb{Z}$Q}$. 
  Consider the isomorphism  $\Phi :\mathcal{Q} \to  \bar{H}^*,e^{\text{n$\alpha $}_1} \mapsto \e^{\text{nQ}}$.  We will identify $\hbs$ with its group algebra $\mathbb{C}[\hbs]$ by $\alpha  \mapsto \e^{\alpha }$.

Just like in the elliptic case considered in \cite{K2}, we identify
$\widehat{f}=f(P) \in\C [H]$ and  meromorphic functions on $\h^*$ by 
\be
%&&\widehat{f}(\mu)=f(<\mu,P>),\quad \mu\in {H}^*
\widehat{f}(\mu)=f<\mu,P>,\quad \mu\in {H}^*
\en
and consider  the field of  meromorphic functions  $M_{{\h}^*}$ on ${\h}^*$ 
\be
&&{M}_{{\h}^*}=\left\{ \widehat{f}:{\h}^*\to \C\ \left|\ 
\widehat{f}=f(P)\in \C [H]\right.\right\}.
\en
\subsubsection{Definition of the $H$-algebra \texorpdfstring{$U_{q,x}(\glth)$}{Uqx}}

Let $\K:=\C[H]$ and define the $H$-algebra $U_{q,x}(\glth):=\K[U_q{(\glth)}] \otimes \C[\bar{H}^*]$.
The moment maps are given by:
%\vspace{-1pt}
\bea
&&\mu_l(\hf)=f(P+h), \qquad   \mu_r(\hf)=f(P). 
\ena
The $H$-bigrading is defined by
%\vspace{-35pt}
\be
&&U_{q,x}(\glth)=\bigoplus_{\alpha,\beta \in H}{U_{q,x}(\glth)_{\alpha,\beta}},\\
&& U_{q,x}(\glth)_{\alpha,\beta}=
\left\{\ x\in U_{q,x}(\glth) \left|\ \mmatrix{q^{P+h}xq^{-(P+h)}=q^{<\alpha, P>}x, \cr
q^{P}xq^{-P}=q^{<\beta, P>}x\cr}\ \right.\right\}.
\en
\noi
For $a,b$ in $\mathbb{C}[U_q(\glth)]$, the multiplication in  $U_{q,x}(\glth)$ is defined through the expression
\bea \lb{mult}
&&\left(f(P)a \otimes  e^{\alpha }\right)\cdot \left(g(P)b \otimes  e^{\beta }\right)=f(P)g(P+<\alpha ,P>)\text{ab}\otimes e^{\alpha +\beta }.
\ena

\begin{dfn}[The Dynamical Currents of $U_{q,x}(\glth)$] \lb{deftotal}
%\vspace{-15pt}
\bea \lb{totalcurrents}
&&E(u)=e(z)e^{2Q}, \quad
F(u)=f(z), \nn\\
&&K_1^\pm(u)=k_1^\pm(z)e^Q, \quad
K_2^\pm(u)=k_2^\pm(z)e^{-Q}, \nn\\
&&H^{\pm }(u)=K_1^{\pm }(u)K_2^{\pm }(u)^{-1}.
%&&H^{\pm }(u)=K_1^{\pm }\left(u-\frac{c}{4}\right)K_2^{\pm }\left(u-\frac{c}{4}\right)^{-1}
\ena
\end{dfn}

\begin{prop}\lb{Defrelns} Let $\eta (a): = 1-q^{2a}$, then the generators of $U_{q,x}(\glth)$ satisfy:
\vspace{-6pt}
\be
&&q^c:\hbox{ {\rm central}}, \lb{u1}\\
&& [h,E(u)]=2E(u),\quad [h,F(u)]={-2}F(u),\\
&&K_1^{\pm }(u)K_1^{\pm }(v)=K_1^{\pm }(v)K_1^{\pm }(u), \lb{u3} \quad
K_2^{\pm }(u)K_2^{\pm }(v)=K_2^{\pm }(v)K_2^{\pm }(u), \lb{u4}
\quad
K_1^{\pm }(u)K_2^{\pm }(v)=K_2^{\pm }(v)K_1^{\pm }(u),
\lb{u2}\\
&&K_1^+(u)K_1^-(v)=\frac{\rho^{+} \left(u-v-\frac{c}{2}\right)}{\rho^{+} \left(u-v+\frac{c}{2}\right)}K_1^-(v)K_1^+(u),\lb{u5}\quad
%&&K_1^-(u)K_1^+(v)=\frac{\rho \left(u-v+\frac{c}{2}\right)}{\rho \left(u-v-\frac{c}{2}\right)}K_1^+(v)K_1^-(u),\lb{u6}\\
K_2^+(u)K_2^-(v)=\frac{\rho^{+} \left(u-v-\frac{c}{2}\right)}{\rho^{+} \left(u-v+\frac{c}{2}\right)}K_2^-(v)K_2^+(u),\lb{u7}\\
%&&K_1^-(u)K_2^+(v)=\frac{\rho \left(u-v+\frac{c}{2}\right)}{\rho \left(u-v-\frac{c}{2}\right)}\frac{\eta \left(u-v-\frac{c}{2}\right)\eta \left(u-v+\frac{c}{2}-1\right)}{\eta \left(u-v+\frac{c}{2}\right)\eta \left(u-v-\frac{c}{2}-1\right)}K_2^+(v)K_1^-(u),\lb{u8}\\
&&K_1^+(u)K_2^-(v)=\frac{\rho^{+} \left(u-v-\frac{c}{2}\right)}{\rho^{+} \left(u-v+\frac{c}{2}\right)}\frac{\eta \left(u-v+\frac{c}{2}\right)\eta \left(u-v-\frac{c}{2}-1\right)}{\eta \left(u-v-\frac{c}{2}\right)\eta \left(u-v+\frac{c}{2}-1\right)}K_2^-(v)K_1^+(u),\lb{u9}\\
&&K_2^+(u)K_1^-(v)=\frac{\rho^{+} \left(u-v-\frac{c}{2}\right)}{\rho^{+} \left(u-v+\frac{c}{2}\right)}\frac{\eta \left(u-v-\frac{c}{2}\right)\eta \left(u-v+\frac{c}{2}+1\right)}{\eta \left(u-v+\frac{c}{2}\right)\eta \left(u-v-\frac{c}{2}+1\right)}K_1^-(v)K_2^+(u),\lb{u99}\\
&&K_1^\pm(u)^{-1}E(v)K_1^\pm(u)=\frac{q\eta \left(u-v\mp\frac{c}{4}-1\right)}{\eta \left(u-v\mp\frac{c}{4}\right)}E(v),\lb{u10}\quad\
%&&K_1^-(u)^{-1}E(v)K_1^-(u)=\frac{q\eta \left(u-v+\frac{c}{4}-1\right)}{\eta \left(u-v+\frac{c}{4}\right)}E(v),\lb{u11}\\
K_1^\pm(u)F(v)K_1^\pm(u)^{-1}=\frac{q\eta \left(u-v\pm\frac{c}{4}-1\right)}{\eta \left(u-v\pm\frac{c}{4}\right)}F(v),\lb{u12}\\%K_1^-(u)F(v)K_1^-(u)^{-1}=\frac{q\eta \left(u-v-\frac{c}{4}-1\right)}{\eta \left(u-v-\frac{c}{4}\right)}F(v),\lb{u13}\\
&&K_2^\pm(u)^{-1}E(v)K_2^\pm(u)=\frac{q^{-1}\eta \left(u-v\mp\frac{c}{4}+1\right)}{\eta \left(u-v\mp\frac{c}{4}\right)}E(v),\lb{u14}\
%&&K_2^-(u)^{-1}E(v)K_2^-(u)=\frac{q^{-1}\eta \left(u-v+\frac{c}{4}+1\right)}{\eta \left(u-v+\frac{c}{4}\right)}E(v),\lb{u15}\\
K_2^\pm(u)F(v)K_2^\pm(u)^{-1}=\frac{q^{-1}\eta \left(u-v\pm\frac{c}{4}+1\right)}{\eta \left(u-v\pm\frac{c}{4}\right)}F(v),\lb{u16}\\
%&&K_2^-(u)F(v)K_2^-(u)^{-1}=\frac{q^{-1}\eta \left(u-v-\frac{c}{4}+1\right)}{\eta \left(u-v-\frac{c}{4}\right)}F(v),\lb{u17}\\
&&E(u)E(v)=\frac{q^{-2}\eta\left(u-v+1\right)}{\eta \left(u-v-1\right)}E(v)E(u) \lb{u18}, \quad
F(u)F(v)=\frac{q^2\eta \left(u-v-1\right)}{\eta\left(u-v+1\right)}F(v)F(u),\lb{u19}\\
&&\left[E(u),F(v)\right]= (q-q^{-1})\left(\delta \left(q^{c}\frac{z}{w}\right)H ^+\left(q^{\frac{-c}{2}}w\right)-\delta \left(q^{-c}\frac{z}{w}\right) H^-\left(q^{\frac{c}{2}}w\right)\right).\lb{u20}
\en
\end{prop}

\noi
{\em Proof.} Straightforward.  The functions $\rho^\pm(u)$ are defined in \eqref{rho} and \eqref{rhominus}.   Using Definition \ref{deftotal}, one can easily verify the defining relations for ${U_{q,x}(\glth)}$ given in the  proposition. 
\qed

 As suggested (without proof) in Section 5 of \cite{K2}, the elliptic half currents survive the degeneration $p \to 0$. For  $U_{q,x}(\glth)$, our definition of the dynamical trignometric (total) currents in Definition \ref{deftotal} and the half-currents in \ref{halfcurrentintegrals} and \ref{halfcurrentseries} is based on this observation.
%\newpage
\subsubsection{Half Currents}
The elements $e^\pm(z)$ and $f^\pm(z)$ are not given explicitly in the definition of $U_q(\glth)$ in \cite {DF}.  In the dynamical case, they can be conveniently expressed with the aid of the Heisenberg algebra, either as contour integrals of the total currents \ref{deftotal} or as Laurent series in $z$ with a modification in the zero Fourier modes $E_0$ and $F_0$.  We will now give their definitions along with their RLL-type commutation relations.
%\newpage
\begin{dfn}{\bf Positive Half-Currents as Integrals}\lb{halfcurrentintegrals}
\be
&&E^+(u)
=q^{-1}a_1 \oint_{C_1} E(u') 
\frac{q^{2(P-1)}\eta\left(u-u'-c/4-P+1\right)\eta(1)}
{\eta(u-u'-c/4)\eta(P-1)}
\frac{dz'}{2\pi i z'},
\lb{Eplus}\\
&&F^+(u)
=q^{-1}a_2 \oint_{C_2} F(u') 
\frac{\eta(u-u'+c/4+P+h-1)\eta(1)}{\eta(u-u'+c/4)\eta(P+h-1)}
\frac{dz'}{2\pi i z'}. 
\lb{Fplus}
\en
Here the contours are chosen as
\be
&&C_1 : 0<|z'|<|q^{c/2}z|, \qquad
C_2 : 0<|z'|<|q^{-c/2}z|,
\lb{C}
\en
%\vspace{-10pt}
and the constants $a_1,a_2$ are chosen to satisfy 
$q^{-2}a_1a_2\eta (1)^2 = -1.$
%check matching in defn to match in e+f+  half curr proof and RLL main thm last no k factor in H=k K1K2inv for me.

\end{dfn}
\noi
For our purposes, it will sometimes be more convenient to do calculations with the series versions of these definitions.   Having access to the individual modes gives us more flexibility and makes the dynamical structure more transparent.  Note that whereas the total currents do not explicitly depend on $P$, the half-currents do.
The key is the Laurent expansion valid in the domain $|z|>1$:
\bea \lb{etaus}
&& \frac{\eta(u+s)}{\eta(u)\eta(s)}=-\sum_{n\in \Z_{\ge 0}}\frac{1}{1-q^{-2s}\delta_{0,n}}z^{-n}.
\ena
A useful consequence that will be used in the proof of the main Theorem \ref{mainthm} is the following relation: %define delta!!!
\bea \lb{etadelta}
&&\frac{\eta(u+s)}{\eta(u)\eta(s)}+\frac{\eta(-u-s)}{\eta(-u)\eta(-s)}=-\delta(q^{2u}),
\ena
with the first summand expanded about $z=\infty$ and the second about $z=0$.
%This expression \ref{Main}  $[E(u),F(v)]$ relation.% RLL proof. and for the decomposition of the representations of the total currents into half-currents.

Using the expansion \eqref{etaus}, Definition  \ref{halfcurrentintegrals} of the positive half-currents is reformulated as infinite series that can be easily extended to the negative counterparts. 
%\newpage
\begin{dfn}{\bf The Series Representation of Positive and Negative Half-Currents.} \lb{halfcurrentseries}
\vspace{-24pt}
\be
&&E^+(u)=
e^{2Q}q^{-1}a_1\eta(1)\left(e_0\frac{1}{1-q^{-2(1-P)}}+ \sum _{n>0} e_n(q^{-\frac{c}{2}}z)^{-n} \right),
\lb{Epluss} \\
&&E^-(u)
=-e^{2Q}q^{-1}a_1\eta(1)\left(e_0\frac{1}{1-q^{2(1-P)}}+ \sum _{n<0} e_n(q^{\frac{c}{2}}z)^{-n}\right),
\lb{Eminuss}\\
%%%
&&F^+(u)=
-q^{-1}a_2\eta(1)\left(f_0\frac{1}{1-q^{-2(P+h-1)}}+ \sum _{n>0} f_n(q^{\frac{c}{2}}z)^{-n} \right),
\lb{Fpluss} \\
&&F^-(u)
=q^{-1}a_2\eta(1)\left(f_0\frac{1}{1-q^{2(P+h-1)}}+ \sum _{n<0} f_n(q^{-\frac{c}{2}}z)^{-n}\right),
\lb{Fminuss}
\en here the positive (resp. negative) currents are expanded around $\infty$ (resp. 0).

\noi
{\em {This yields the following decomposition of the total currents}}
\end{dfn}
\vspace{-6pt}
%We see thae:
\noi
\bea \lb{totaldecomp}
%&&{\text {clearly}}\nn\\
&&\ \ q^{-1} a_1\eta (1)E(u)=E^+\left(u+\frac{c}{4}\right)-E^-\left(u-\frac{c}{4}\right),\nn\\
&&-q^{-1}a_2\eta(1)F(u)=F^+\left(u-\frac{c}{4}\right)-F^-\left(u+\frac{c}{4}\right).\qquad
\ena
\begin{prop}[Commutation relations for Positive Half-Currents] \lb{plusrelns} Let $u=u_1-u_2$.
\bea
&& K_i^+(u_1)K_j^+(u_2)=K_j^+(u_2)K_i^+(u_1), \qquad   (i,j =1,2) \lb{plusrelkk} \\
&&\kpone^{-1}\Ept\kpone=\Ept\frac{q\eta(u-1)}{\eta(u)}+\Ep\frac{q^{-1}\eta(1)\eta(P+u-2)}{\eta(P-2)\eta(u)},
\lb{plusrelek1}\\
&&\kpone\fpt\kpone^{-1}=\fpt \frac{q\eta(u-1)}{\eta(u)} +\fp\frac{q^
{2u-1}\eta(1)\eta(P+h-u)}{\eta(P+h)\eta(u)},
\lb{plusrelfk1}\\
&& K_2^+ (u_1)^{-1 } E^+(u_2)K_2^+ (u_1)= E^+(u_2)\frac{q^{-1}\eta(u+1)}{{\eta }(u)}-E^+ (u_1)\frac{q^{-1}\eta   (1) \eta   (P+u)}{\eta   (P) \eta   (u)} , 
\lb{plusrelek2} \\
%&& K_2^+ (u_1)F^+ (u_2)  K_2^+ (u_1)^{-1}= F^+(u_2)\frac{\eta(u+1)}{{q\eta }(u)}-F^+ (u_1)\frac{\eta   (1) \eta   (P+u)}{q\eta   (P) \eta   (u)} , 
&&K_2^+(u_1)F^+(u_2)K_2^+(u_1){}^{-1}= F^+(u_2)\frac{q^{-1}\eta (u+1)}{\text{$\eta $}(u)} -\fp\frac{q^{2u-1}\eta (1)\eta (P+h-u-2)}{\eta (P+h-2)\eta (u)},
\lb{plusrelfk2} \\
&&\frac{q^{2u-1}\eta (1-u)}{\eta (u)}\hspace{0.2cm}E^+\left(u_1\right)E^+\left(u_2\right)+\frac{q^{-1}\eta (1+u)}{\text{$\eta $}(u)}\hspace{0.2cm}E^+\left(u_2\right)E^+\left(u_1\right) \nn\\
&&\hspace{0.5cm} =E^+\left(u_1\right){}^2\hspace{0.2cm}\frac{q^{-1}\eta (1)\eta (P+u-2)}{\text{$\eta $}(u)\eta (P-2)}+E^+\left(u_2\right){}^2\hspace{0.2cm}\frac{q^{2u-1}\eta (1)\eta (P-u-2)}{\text{$\eta $}(u)\eta (P-2)},
\lb{plusrelee} \\
%\ena
%\bea
&&\frac{q^{-1}\eta (1+u)}{\eta (u)}\hspace{0.2cm} F^+\left(u_1\right) F^+\left(u_2\right)+\frac{q^{2u-1}\eta (1-u)}{\eta (u)}\hspace{0.2cm} F^+\left(u_2\right) F^+\left(u_1\right)  \nn\\
&&\hspace{0.5cm}=F^+\left(u_1\right){}^2\hspace{0.2cm}\frac{q^{2u-1}\eta (1)\eta (P+h-u-2)}{\eta (P+h-2)\eta (u)}+ F^+\left(u_2\right){}^2\hspace{0.2cm}\frac{q^{-1}\eta (1)\eta (P+h+u-2)}{\text{$\eta $}(P+h-2)\eta (u)},
 \lb{plusrelff}\\
%&&{and} \nn\\
&&\left[E^+\left(u_1\right),F^+\left(u_2\right)\right] \nn\\
&&\quad=q^{2u}\left(q^{-1}-q\right)\left(K_2^+\left(u_2\right){}^{-1}K_1^+\left(u_2\right)\frac{\eta (P-u-1)}{\eta (u)\eta (P-1)} -K_2^+\left(u_1\right){}^{-1}K_1^+\left(u_1\right)\frac{\eta (P+h-u-1)}{\eta (u)\eta (P+h-1)}\right).
\lb{plusrelef} %won't fit!
\nn\\
\ena
\end{prop}
\noi
{\it Proof.}  
%%%%%%
 The relations \eqref{plusrelkk} are direct consequences of the definitions of $K_i^+(u)$, for i=$1,2$. 
 The remaining equations can be grouped into three types :

%\begin{enumerate}
\noi
(i) {\eqref{plusrelek1} - \eqref{plusrelfk2}}.
By using the defining relations in Proposition \ref{Defrelns} along with the contour integral definitions of the half-currents \ref{halfcurrentintegrals} and the following identity, one can verify the commutation relations for $E^+(u_2)$ or $F^+(u_2)$ with $K_1^{+}(u_1)$ and  $K_2^{+}(u_1):$% in equations \eqref{plusrelek2} - \eqref{plusrelfk1}.  
\bea \lb {etaid3}
&&\frac{\eta(u_1+t)\eta(u_2 + s)}{\eta(u_1) \eta(u_2) \eta(s)} =  \frac{\eta(u_1- u_2+t)\eta(u_2 + s + t)}{\eta(u_1 - u_2) \eta(u_2) \eta(s + t)}+
\frac{\eta(u_2- u_1 + s) \eta(u_1 + s + t) \eta(t)}{\eta(u_2-u_1)\eta(u_1) \eta(s) \eta(s + t)}.\nn\\
\ena
(ii) {Relations \eqref{plusrelee} and \eqref{plusrelff}}.   To prove \eqref{plusrelee}, symmetrize the integration variables by writing each term in \eqref{plusrelee} as a sum of two integrals (each one in both $u_1'$ and $u_2'$), with the $E$-part $E(u_1')E(u_2')$ (using Proposition \ref{Defrelns}), combine the first term on the LHS and first term on the RHS of \eqref{plusrelee}, and use the fact (eqn.(4.16) in \cite{JKOS2}) that the following expression remains unchanged on interchanging $u_1$ and $u_2$ only:

\be
&& -  \frac{q^{-3t}\eta(t)\eta(u_2-u_1+s)\eta(u_1-u_1'+s+t)\eta(t)\eta(u_1-u_2'+s-t)\eta(t)}{\eta(s)\eta(s-t)\eta(s+t)\eta(u_2-u_1)\eta(u_1-u_1')\eta(u_1-u_2')} \nn\\
&& + \frac{ q^{(-3t+u_2-u_2')}\eta(t)\eta(u_2'-u_1'+t)\eta(u_1-u_2'+s+t)\eta(u_1-u_2-t)\eta(u_2-u_1'+s-t)\eta(t)}{\eta(s-t)\eta(s+t)\eta(u_1-u_2)\eta(u_1-u_2')\eta(u_2-u_1')\eta(u_2-u_1'-t)}\nn\\
&& -  \frac{q^{(-5t-u_1-u_1')}\eta(t)\eta(u_2-u_1+s)\eta(u_1-u_1'+s-t)\eta(t)\eta(u_1-u_2'+s+t)\eta(u_2'-u_1'+t)\eta(t)}{\eta(s)\eta(s-t)\eta(s+t)\eta(u_2-u_1)\eta(u_1-u_1')\eta(u_1-u_2')\eta(u_2'-u_1-t)}\nn\\
&&+\frac{q^{-t}\eta(t)\eta(u_1-u_2-t)\eta(u_1-u_1'+s+t)\eta(u_2-u_2'+s-t)\eta(t)}{\eta(s-t)\eta(s+t)\eta(u_1-u_2)\eta(u_1-u_1')\eta(u_2-u_2')},
\en
setting $s=-P+2, t=1, u_1=u_1-\frac{c}{4}, u_2=u_2-\frac{c}{4}$.  The proof of \eqref{plusrelff} is similar.

\noi
(iii) Finally, use Definition  \ref{halfcurrentseries} for the positive half-currents series to obtain the $[E^+(u_1),F^+(u_2)]$ relation \eqref{plusrelef}.   Since $\left[E^+(u),q^{-2(P+h-1)}\right]=0$, we have
\bea \lb{efpluscheck}
&& \left[E^+(u_1),F^+(u_2)\right]= -{e^{2Q}q^{-2}a_1a_2\eta (1)^2} \Biggl( \ \left[e_0\frac{q^{2(P-1)}}{1-q^{2(P-1)}}+\sum _{n=0}^{\infty } e_n(q^{-\frac{c}{2}}z_1)^{-n}, f_0\right]\frac{1}{1-q^{-2(P+h-1)}}  \nn\\
&&\qquad  \qquad \qquad \qquad \ +\left[e_0\frac{q^{2(P-1)}}{1-q^{2(P-1)}}+\sum _{n=0}^{\infty } e_n(q^{-\frac{c}{2}}z_1)^{-n}                                                       , \sum _{n=1}^{\infty } f_n(q^{\frac{c}{2}}z_2)^{-n}\right] \ \Biggr).\nn\\
\ena

%\vspace{-32pt}
\noi
 Expand the four terms on the right side of \eqref{efpluscheck}, using the commutation relation between $e_n$ and $f_n$ in the defining relation
\eqref{efcomm}:%check matching in defn to match in e+f+  half curr proof and RLL main thm last

\be
&&\left[\sum _{n=0}^{\infty } e_n(q^{-\frac{c}{2}}z_1)^{-n},\sum _{n=1}^{\infty } f_n(q^{\frac{c}{2}}z_2)^{-n}\right]=\frac{q-q^{-1}}{1-q^{2(u_2-u_1)}}\left(\psi \left(z_2\right)-\psi \left(z_1\right)\right), \\
&&\left[e_0\frac{q^{2(P-1)}}{1-q^{2(P-1)}},f_0\right]\frac{1}{1-q^{-2(P+h-1)}}= \frac{(q-q^{-1})q^{2(P-1)}}{(1-q^{2(P-1)})(1-q^{-2(P+h-1)})}\left(q^h-q^{-h}\right),\\
\en
\be
&&\left[\sum _{n=0}^{\infty } e_n(q^{-\frac{c}{2}}z_1)^{-n}, f_0\right]\frac{1}{1-q^{-2(P+h-1)}}=\frac{q-q^{-1}}{1-q^{-2(P+h-1)}}\left(\psi \left(z_1\right)-q^{-h}\right),\\
&&\left[e_0\frac{q^{2(P-1)}}{1-q^{2(P-1)}},\sum _{n=1}^{\infty } f_n(q^{\frac{c}{2}}z_2)^{-n}\right]=\frac{(q-q^{-1})q^{2(P-1)}}{1-q^{2(P-1)}}\left(\psi \left(z_2\right)-q^h\right).
\en
Now use these four expressions in  Equation \eqref{efpluscheck}, along with the formula:
\bea
&& \frac{1}{\eta(x)}-\frac{1}{\eta(y)}=\frac{\eta(x-y)}{\eta(x)\eta(-y)}.
\ena
\qed\\
\noi
{\em Remark.} We do not solve the commutation relations for the {\em negative} half-currents because we will not explicitly need them in the sequel.

\section{Definition of \texorpdfstring{$U(R_{q,x}(\glth))$}{U(R)} and the Main Theorem}
\subsection{Definition of $U(R)$}
We will use the following presentation of the dynamical affine $R$-matrix, obtained as the degeneration of the elliptic $R$-matrix in \eqref{Rmat4} as $p \to 0$ (see \cite{idzumi}, \cite{K2}, also eq(5.9) and the remark below it, in \cite{clad}) :
\bea \lb{rplus}
&&R^\pm(u,P)=\rho^\pm (z)\left(
\begin{array}{cccc}
 1 & 0 & 0 & 0 \\
 0 & \frac{\text{q$\eta $}(P+1)\eta (P-1)\eta (u)}{\eta (P)^2\eta (u+1)} & \frac{\eta (1)\eta (P+u)}{\eta (u+1)\eta (P)} & 0 \\
 0 & \frac{\eta (1)\eta (P-u)q^{2u}}{\eta (u+1)\eta (P)} & \frac{\text{q$\eta $}(u)}{\eta (u+1)} & 0 \\
 0 & 0 & 0 & 1 \\
\end{array}
\right) \nn\\
&&=\rho^\pm (z)\left(
\begin{array}{cccc}
 1 & 0 & 0 & 0 \\
 0 & \frac{q\left(1-q^{2(P-1)}\right) \left(1-q^{2(P+1)}\right)\left(1-q^{2u}\right)}{\left(1-q^{2P}\right)^2 \left(1-q^{2(u+1)}\right)} & \frac{\left(1-q^2\right) \left(1-q^{2(P+u)}\right)}{\left(1-q^{2P}\right) \left(1-q^{2(u+1)}\right)} & 0 \\
 0 & \frac{q^{2u}\left(1-q^2\right)\left(1-q^{2(P-u)}\right)}{\left(1-q^{2 P}\right)\left(1-q^{2(u+1)}\right)} & \frac{q\left(1-q^{2u}\right)}{1-q^{2(u+1)}} & 0 \\
 0 & 0 & 0 & 1 \\
\end{array}
\right).
\ena
Here, the coefficient is given by
\bea \lb{rho}
&&\rho^+ (z)=\frac{\sqrt{q} \left(z^{-1};q^4\right)_\infty \left(q^4z^{-1};q^4\right)_\infty}{\left(q^2z^{-1};q^4\right)_\infty^2},\quad \text{}
 (a;x)_\infty=\prod_{k=0}^\infty(1-ax^{k})
%&&\Theta_{p}(z)=(z;p)_\infty(p/z;p)_\infty(p;p)_\infty,
\ena
and $\rho^-(z)=(\rho^+(z^{-1}))^{-1}$, which expands as
\bea \lb{rhominus}
&&\rho^- (z)=\frac{q^{-\frac{1}{2}}{\left(q^2 z;q^4\right)_\infty^2}}{ \left(z;q^4\right)_\infty \left(q^4z;q^4\right)_\infty}.
%&&\Theta_{p}(z)=(z;p)_\infty(p/z;p)_\infty(p;p)_\infty,
\ena
%and satisfies the consistency condition: $\rho (z)\rho \left(z^{-1}\right)=1$. 
We will write $\rho^\pm(u)$ for $\rho^\pm(z)=\rho^\pm(q^{2u})$ when the context is clear.  We then have, $\rho^+ (u)\rho^- (-u)=1.$

  Observe that the quantum affine $R$-matrix  $R_q(z)$ in Section \ref{subsecdf} is the degeneration of $R^+(u,P)$ as $P \to \infty$.{\footnote {The factors $\rho^\pm(z)$ can be chosen arbitrarily and are set to 1 in $R_q(z)$.} }% and that Ding-Frenkel do not use any $\rho(z)$-type coefficient in their formulation. \\{\text{and}}
It is easily verified that the {\em unitarity} condition (eq(4.9) in \cite{DF}) continues to hold in the dynamical case:
\vspace{-.3cm}
\be
&& R^+_{21}(u,P)^{-1 }=R^-(-u,P).
\en
where the inverted $R$-matrix is
\bea R^+(u,P)^{-1}=\rho^+(z)^{-1}
\left(
\begin{array}{cccc}
 1 & 0 & 0 & 0 \\
 0 & \frac{\eta (u)}{q \eta (u-1)} & -\frac{\eta (1) \eta (P+u)}{q^2 \eta (u-1) \eta (P)} & 0 \\
 0 & \frac{\eta (1) \eta (P-u)}{\eta (1-u) \eta (P)} & \frac{\eta (u) \eta (P-1) \eta (P+1)}{q \eta (u-1) \eta (P)^2} & 0 \\
 0 & 0 & 0 & 1 \\
\end{array}
\right).
\ena 
As mentioned in the introduction, Arnaoudon et al.\cite{clad} suggest a definition of $U_{q,\la}(\slth)$ based on $RLL$ relations. Besides being motivated by their considerations, the definition presented in this treatise is natural, from the location of $R_{q,x}$ (resp. $U_{q,x}(\slth)$) in the cladistics of quantum $R$-matrices (resp. quantum algebras) as illustrated in Figure \ref{fig:figure1})  (resp. Figure \ref{fig:algebradegenerations}). 
We can extend the Heisenberg algebra for $U_{q,x}(\slth)$ to include the element $h$ as follows:
\begin{dfn} {Heisenberg algebra for $U(R)$}.  Let
\bea
&&H=\C P + \C (P+h), \quad H_Q^*=\sum \C Q, \quad H^*=\frak{h}^* + H_Q^* \nn
\ena
\noi
We identify $\widehat{f}=f(P,P+h) \in\C [H]$ and  meromorphic functions on $\h^*$
via
\bea
&&f(P,P+h)(\xi)=f(<P,\xi>,<P+h,\xi>),\quad \xi \in H^*
\ena
\end{dfn}

\begin{dfn} \lb{UR}
 $U(R):=U(R_{q,x}(\glth))$ is the algebra over $M_{H^*}$ with generators $L^\pm(u)$ given by $L^{\pm}(u)=\left(L^\pm_{ab}(u)\right)_{a,b=1}^2$, with $L^\pm_{ab}(u)=\sum_{n=0}^\infty L_{ab,\pm n}q^{\mp 2un}$% such that $L_{aa,o}^+ L_{aa,o}^- = L_{aa,o}^- L_{aa,o}^+ =1$ 
and the RLL-relations:
\bea
&&R^{\pm(12)}(u,P+h)L^{\pm (1)}\left(u_1\right)L^{\pm (2)}\left(u_2\right)=L^{\pm (2)}\left(u_2\right)L^{\pm (1)}\left(u_1\right)R^{\pm(12)}(u,P),  \lb{rll1}\\
&&R^{\pm(12)}\left(u\pm\frac{c}{2},P+h\right)L^{\pm(1)}\left(u_1\right)L^{\mp(2)}\left(u_2\right)=L^{\mp(2)}\left(u_2\right)L^{\pm(1)}\left(u_1\right)R^{\pm(12)}\left(u\mp\frac{c}{2},P\right) \lb{rll2}\nn\\
\ena
\noi
with $u=u_1-u_2$.
\end{dfn}
\begin{lem} \lb{halgUR}
$U(R)$ is an $H$-algebra with the $H$-bigrading defined by
\be
&& U(R)=\bigoplus_{\alpha,\beta \in H}{U(R)_{\alpha,\beta}}\hspace{.1cm},
\quad
 U(R)_{\alpha,\beta}=
\left\{\ x\in U(R) \left|\ \mmatrix{q^{P+h}xq^{-(P+h)}=q^{<\alpha, P>}x, \cr
q^{P}xq^{-P}=q^{<\beta, P>}x\cr}\ \right.\right\},
\en
and moment maps 
\bea
&&\mu_l(f)=f(P+h), \qquad   \mu_r(f)=f(P), \qquad f \in M_{H^*} 
\ena
\end{lem}

\noi 
{\em {Proof}.} A straightforward verification.

\noi
{\em {Remark.}} The action on the generators
 is given by
\bea
f(P+h)L^{\pm}_{ab}(u)&=&
L^{\pm}_{ab}(u)f(P+h-w(a)),\nn\\
f(P)L^{\pm}_{ab}(u)&=&L^{\pm}_{ab}(u)
f(P-w(b)),\lb{shifts}
\ena
where the weight function $w:\left\{1,2\right\} \rightarrow \left\{\pm1\right\}$ is given by identifying 
\bea
 L^{\pm}(u) \to \left(
\begin{array}{cc}
 L_{++}(u) & L_{+-}(u) \\
L_{-+}(u) & L_{--}(u) \\
\end{array}
\right)
\ena
\noi
We also define the subalgebras $U(R^{\pm}):= U\left(R^\pm_{q,x}(\glth)\right)$.
The expansion of the $RLL$-relations \eqref{rll1} and \eqref{rll2} is given in Appendix A.
%%%%%
%%%main thm%%%
%%%%%%

For proving our main result it will be convenient to summarize the left action of $P$ and $P+h$ on the half-currents:
\vspace{-.4cm}
\begin{center}
    \begin{tabular}{ | l | l | l | l | l | l | l | l |}
    \hline
   &$e^Q$ &  $E(u),E^\pm(u)$ & $F(u), F^\pm(u)$ & $K_1^\pm(u)$ & $K_2^\pm(u)$ & $H^+(z)$ & $H^-(z)$ \\ \hline
    $P$ & $P-1$ & $P-2$ & $P$ & $P-1$ & $P+1$ & $P-2$ & $P-2$\\ \hline
   $P+h$ & $P+h-1$ & $P+h$ & $P+h-2$ & $P+h-1$ & $P+h+1$ & $P+h-2$ & $P+h-2$\\ \hline
     \end{tabular}
\end{center}
For example, the third column can be interpreted as the following relations:

\noi
 $f(P)E(u)=E(u)f(P-2)$, $f(P+h)E(u)=E(u)f(P+h)$,  
$f(P)E^\pm(u)=E^\pm(u)f(P-2)$ and $f(P+h)E^\pm(u)=E^\pm(u)f(P+h)$.
%\newpage
\subsection{ Isomorphism of $H$-Algebras \texorpdfstring{$ U_{q,x}(\glth) \simeq U({R_{q,x}}(\glth))$.}{Uqx and U(R)}} 
We now prove our main result, an extension of the Ding-Frenkel type isomorphism (see Subsection \ref{subsecdf}) to the dynamical case.
\begin{thm}[Main Theorem] \lb{mainthm}
There exists a unique Gauss decomposition of $L^\pm(u)$:
%\vspace{-12pt}
\bea \label{Lplusop}
&&L^\pm(u)=\left(
\begin{array}{cc}
 1 & F^\pm(u) \\
 0 & 1 \\
\end{array}
\right)\left(
\begin{array}{cc}
 K_1^\pm(u) & 0 \\
 0 & K_2^\pm(u) \\
\end{array}
\right)\left(
\begin{array}{cc}
 1 & 0 \\
 E^\pm(u) & 1 \\
\end{array}
\right) \\
&&\text{  }\text{  }\text{  }\text{  }\text{  } \text{  }\text{  } \text{  }\text{  }  =\left(
\begin{array}{cc}
 K_1^\pm(u)+F^\pm(u)K_2^\pm(u)E^\pm(u) & F^\pm(u)K_2^\pm(u) \\
 K_2^\pm(u)E^\pm(u) & K_2^\pm(u) \\
\end{array}
\right), \lb{gaussL}
\ena
%\vspace{-48pt}

\noi
which yields an an isomorphism of $M_{H^*}$-algebras:
%\vspace{-36pt}
\be
&&\Phi : U_{q,x}\left(\glth\right) \longrightarrow U\left(R_{q,x}(\glth)\right), \text{ defined by}\\
&&\hspace{36pt} E(u) \mapsto  E^+\left(u+\frac{c}{4}\right)-E^-\left(u-\frac{c}{4}\right),\\
&&\hspace{36pt}F(u) \mapsto  F^+\left(u-\frac{c}{4}\right)-F^-\left(u+\frac{c}{4}\right),\\
&&\hspace{30pt}K_i^{\pm }(u) \mapsto  K_i^{\pm }(u), \quad q^{\pm c }\mapsto  q^{\pm c }.
\en
\end{thm}
%\vspace{-48pt}
\noi
{\it Proof.} 
To prove that $\Phi$ is a homomorphism of associative algebras, we will need the inverted $L$-operators, which are obtained using the Gauss decomposition, as:
\bea \lb{linv}
\left(L^\pm(u)\right)^{-1}= \left(
\begin{array}{cc}
 K_1^\pm(u)^{-1} & -K_1^\pm(u)^{-1}F^\pm(u) \\
 -E^\pm(u)K_1^\pm(u)^{-1} & K_2^\pm(u)^{-1}+E^\pm(u)K_1^\pm(u)^{-1}F^\pm(u) \\
% K^\pm(u-1)^{-1} & -K^\pm(u-1)^{-1}F^\pm(u) \\
 %-E^\pm(u)K^\pm(u-1)^{-1} & K^\pm(u)+E^\pm(u)K^\pm(u-1)^{-1}F^\pm(u) \\
\end{array}
\right).
\ena

%%%%%
%%%PROOF OF  MAIN THM
\noi
The verification reduces to checking that the $RLL$ relations \eqref{rll1} and \eqref{rll2} imply the defining relations \eqref{Defrelns}.  Two key differences from Theorem \ref{dfthm}  are that $K_i^{+ }(u)$ and  $K_i^{-}(u)$ do not commute, for $i=1$ and 2, and that the entries in $L^\pm(u)$ (the half-currents) 
do not commute with the entries of the $R$-matrix.  The $RLL$-relations \eqref{rll1} and \eqref{rll2} are inverted to yield:

\bea \lb{invertedrll}
&&L^{\pm (1)}(u_1){}^{-1 }R^{\pm(12)}(u,P+h)^{-1}L^{\pm (2)}(u_2) = L^{\pm (2)}(u_2)R^{\pm(12)}(u,P)^{-1}L^{\pm (1)}(u_1){}^{-1}, \lb{eq4.2} \\
&&L^{-(1)}(u_1){}^{-1 }R^{-(12)}(u-\frac{c}{2},P+h)^{-1}L^{+(2)}(u_2) = L^{+(2)}(u_2)R^{-(12)}(u+\frac{c}{2},P)^{-1}L^{-(1)}(u_1){}^{-1}, \lb{eq4.3} \nn\\ \\
&&L^{+(1)}(u_1){}^{-1 }R^{+(12)}(u+\frac{c}{2},P+h)^{-1}L^{-(2)}(u_2) = L^{-(2)}(u_2)R^{+(12)}(u-\frac{c}{2},P)^{-1}L^{+(1)}(u_1){}^{-1}, \lb{eq4.5} \nn\\ \\
&&R^{+(12)}(u+\frac{c}{2},P+h)^{-1}L^{-(2)}(u_2)L^{+(1)}(u_1) = L^{+(1)}(u_1)L^{-(2)}(u_2)R^{+(12)}(u-\frac{c}{2},P)^{-1}\text{  }, \lb{eq4.4} \\
%\ena
%\bea
&&L^{\pm (2)}(u_2){}^{-1}L^{\pm (1)}(u_1){}^{-1}R^{\pm(12)}(u,P+h)^{-1} = R^{\pm(12)}(u,P)^{-1}L^{\pm (1)}(u_1){}^{-1}L^{\pm (2)}(u_2){}^{-1}, \lb{eq4.6} \\
&&L^{+(2)}(u_2){}^{-1}L^{-(1)}(u_1){}^{-1}R^{-(12)}(u-\frac{c}{2},P+h)^{-1} = R^{-(12)}(u+\frac{c}{2},P)^{-1}L^{-(1)}(u_1){}^{-1}L^{+(2)}(u_2){}^{-1}.  \lb{eq4.7} \nn\\ 
\ena
%The expansion of these relations is illustrated in Appendix B.
(i).  Let $(A)[[i,j]]$ denote the elements in row $i$ and column $j$ of the matrix equation labelled as $(A)$.  The relations involving only the $K_i^\pm(u)$  are easily obtained by directly equating the matrix entries in both  sides of \eqref{eq4.7}[[4,4]], \eqref{rll2}[[4,4]], \eqref{eq4.3}[[2,2]], \eqref{eq4.5}[[2,2]] and using the definition of $\rho^\pm(z)$ in \eqref{rho} and \eqref{rhominus}.

\noi
(ii). Let us prove the relation between $K_i^+(u_1)$ and $E(u_2)$, for $i=1,2$. Use the relation between $\kpone$ and $\kmtwot$ in $\eqref{eq4.5}[[2,1]]$:  \footnote{The   positioning of the $R$-matrix elements is important because they do not commute with the half currents.} 
%for k1 with eminus use 45rho green ya %
\be %45rho
&&\kpone^{-1}\frac{q^{-1}\eta(u+\frac{c}{2}) \rho^+(u+\frac{c}{2})^{-1}}{\eta(u+\frac{c}{2}-1)} \kmtwot\Emt\\
&&\qquad =\kmtwot\frac{q^{-2}\eta(1)\eta(P+u-\frac{c}{2})\rho^+(u-\frac{c}{2})^{-1}}{\eta(P)\eta(u-\frac{c}{2}-1)}\Ep\kpone^{-1} \\
&&\qquad \qquad +\kmtwot\Emt \rho^+(u-\frac{c}{2})^{-1}\kpone^{-1},
\en
which implies that
\be
&&\kpone^{-1}\Emt\kpone=\frac{q^{-1}\eta(1)\eta(P+u-\frac{c}{2})}{\eta(P)\eta(u-\frac{c}{2})}\Ep+\frac{q\eta(u-\frac{c}{2}-1)}{\eta(u-\frac{c}{2})}\Emt.
\en
Similarly, from \eqref{eq4.6}[[2,1]], it follows that
\be %46pp
&&-\Ept\kponet^{-1}\kpone^{-1}\\
&&\qquad =\frac{q^{-2}\eta(1)\eta(P+u)}{\eta(P)\eta(u-1)}\Ep\kpone^{-1}\kponet^{-1}-
\frac{q^{-1}\eta(u)}{\eta(u-1)}\kpone^{-1}\Ept\kponet^{-1},
\en
from which we get
\be
&&\kpone^{-1}\Ept\kpone=\frac{q^{-1}\eta(1)\eta{(P+u)}}{\eta(P)\eta(u)}\Ep +\frac{q\eta(u-1)}{\eta(u)}\Ept.
\en
Hence we conclude that 
\be
&&K_1^+(u_1)^{-1}E(u_2)K_1^+(u_1)=\frac{q\eta \left(u-\frac{c}{4}-1\right)}{\eta \left(u-\frac{c}{4}\right)}E(u_2),
\en
as required.
 Start with the following expression which we obtain from\eqref{rll1}[[4,3]]:
\be
&&K_2^+ (u_1)  K_2^+(u_2)  E^+(u_2) \nn\\
&& =K_2^+(u_2)  E^+ (u_2)   K_2^+ (u_1)  \frac{q \eta   (u )}{\eta   (u+1 )}+ K_2^+  (u_2)   K_2^+ (u_1)  E^+ (u_1)  \frac{\eta   (1) \eta   (P+u)}{\eta (P) \eta   (u+1)} \nn\\
&&\implies K_2^+ (u_1)^{-1}E^+ (u_2)  K_2^+ (u_1)=-E^+ (u_1)\frac{\eta   (1) \eta   (P+u)}{q\eta   (P) \eta   (u)} + E^+(u_2)\frac{\eta(u+1)}{{q\eta }(u)}.
\en
The corresponding entry in  the second equation from \eqref{rll2}[[4,3]] reads:
\be %rllpm so rhoplus ok
&&\rho^+(u+\frac{c}{2})  K_2^+(u_1)   K_2^-(u_2)  E^-(u_2)=K_2^- (u_2)  E^-(u_2)  K_2^+(u_1)  \frac{q \eta(u-\frac{c}{2}) \rho^+ (u-\frac{c}{2})}{\eta(u-\frac{c}{2}+1)}+ \nn\\
&&\qquad K_2^-(u_2)  K_2^+(u_1)  E^+(u_1)  \eta(1)\eta(u-\frac{c}{2}+P) \frac{\rho^+(u-\frac{c}{2})}{\eta(P)}\nn\\
&&\implies K_2^+(u_1)^{-1}E^-(u_2)  K_2^+(u_1)=-E^+(u_1)\frac{\eta (1) \eta \left(P+u-\frac{c}{2}\right)}{\text{q$\eta $}(P) \eta \left(u-\frac{c}{2}\right)} + E^-(u_2)\frac{\eta \left(u-\frac{c}{2}+1\right)}{ \text{q$\eta $}\left(u-\frac{c}{2}\right)}.
\en
Hence, we conclude that 
\be
&&K_2^+(u_1)^{-1}E(u_2)K_2^+(u_1)=q^{-1}\frac{\eta \left(u-\frac{c}{4}+1\right)}{\eta \left(u-\frac{c}{4}\right)}E(u_2),
\en
as required. 

%%%%%%%%%%%%%%%%%%%%%%%%%%%%%%%%
Since the $RLL$ relations are symmetric in $c$, the relations for $K_i^-(u_1)$ with $E(u_2)$ can be obtained by choosing the "negative" counterparts of the expressions in the preceding proof for $K_i^+(u_1)$.  The matrix entries are given in \eqref{eq4.6}[[4,3]] and \eqref{eq4.7}[[4,3]] for $K_1^-$ and in \eqref{rll1}[[4,3]] and \eqref{rll2}[[4,3]] for $K_2^-$.  Note that this reduces to just replacing $c$ by $-c$ in the preceding proof.  Similarly, for the proof of the relation involving $F(u_1)$ and $K_2^\pm(u_2)$, use  the following expressions, obtained from \eqref{rll1}[[3,4]] and \eqref{rll2}[[3,4]]:
\bea \lb{commk2fpm}
&&\kpmtwo\fpmt\kpmtwo^{-1}=-\frac{q^{2u-1}\eta(1)\eta(P+h-u)}{\eta(P+h)\eta(u)}\fpm+\frac{q^{-1}\eta(u+1)}{\eta(u)}\fpmt, \nn\\
&&\kpmtwo\fmpt\kpmtwo^{-1}=-\frac{q^{2u\pm c-1}\eta(1)\eta(P+h-u\mp\frac{c}{2})}{\eta(P+h)\eta(u\pm\frac{c}{2})}\fpm+\frac{q^{-1}\eta(u\pm\frac{c}{2}+1)}{\eta(u\pm\frac{c}{2})}\fmpt \nn\\
\ena
to obtain the required relation.  The verification of the relation between $F(u)$ with $K_1^\pm(u)$ is similar, using \eqref{eq4.5}[[2,1]], \eqref{eq4.6}[[2,1]], and \eqref{eq4.7}[[2,1]]. 

(iii).
We prove the relation between $F(u_1)$ and $F(u_2)$.  The proof for $E(u_1)$ with $E(u_2)$ is similar. 
%%ff%%
Start with the 2 pairs of equalities from \eqref{rll1}[[1,4]] and  \eqref{rll2}[[1,4]] :
\bea
&&\fpm\kpmtwo\fpmt\kpmtwot=\fpmt\kpmtwot\fpm\kpmtwo, \nn\\
&&\rho^\mp(u\pm c/2)\fpmt\kpmtwot\fmp\kmptwo=\fmp\kmptwo\fpmt\kpmtwot\rho^\mp(u\mp c/2)\nn\\.
%&&F^{\pm}(u_2)k_2^{\pm}(u_1)
\ena
Simplifying these using the four expressions \eqref{commk2fpm} yields 4 equations :
\be
&&\frac{q^{-1}\eta (1+u)}{\eta (u)}\hspace{.2cm} F^\pm(u_1) F^\pm(u_2)+\frac{q^{2u-1}\eta (1-u)}{\eta (u)}\hspace{.2cm} F^\pm (u_2) F^\pm (u_1)  \nn\\
&&\hspace{2cm}=F^\pm (u_1)^2\hspace{.2cm}\frac{q^{2u-1}\eta (1)\eta (P+h-u-2)}{\eta (P+h-2)\eta (u)}+ F^\pm (u_2)^2\hspace{.2cm}\frac{q^{-1}\eta (1)\eta (P+h+u-2)}{\text{$\eta $}(P+h-2)\eta (u)}, 
\en
\be
&&\frac{q^{-1}{\rho(u\mp\frac{c}{2})}\eta (u\mp\frac{c}{2}+1)}{{\rho(u\pm\frac{c}{2})}\eta (u\mp\frac{c}{2})}\hspace{.2cm} F^\mp (u_1) F^\pm (u_2)-\frac{q^{-1}\eta (1-u\pm\frac{c}{2})}{\eta (-u\pm\frac{c}{2})}\hspace{.2cm} F^\pm (u_2) F^\mp (u_1)  \nn\\
&&\hspace{2cm}=F^\mp (u_1)^2\hspace{.2cm}\frac{\rho^+(u\mp\frac{c}{2})q^{2u\mp c-1}\eta (1)\eta (P+h-u\pm\frac{c}{2}-2)}{\rho^+(u\pm\frac{c}{2})\eta (P+h-2)\eta (u\pm\frac{c}{2})}\\
&&\hspace{4cm}- F^\pm (u_2)^2\hspace{.2cm}\frac{q^{-2u\mp c-1}\eta (1)\eta (P+h+u\mp\frac{c}{2}-2)}{\text{$\eta $}(P+h-2)\eta (u\mp\frac{c}{2})}.
\en
%%%%%%%%%%%%%%
Now apply these four relations in the expansion of
\be
\eta(u+1)F(u_1)F(u_2)=\eta(u+1)\Bigg(F^+\left(u_1-\frac{c}{4}\right)-F^-\left(u_1+\frac{c}{4}\right)\Bigg)\Bigg( F^+\left(u_2-\frac{c}{4}\right)-F^-\left(u_2+\frac{c}{4}\right)\Bigg),
\en
\be
q^2\eta(u-1)F(u_2)F(u_1)=q^2\eta(u-1)\Bigg(F^+\left(u_2-\frac{c}{4}\right)-F^-\left(u_2+\frac{c}{4}\right)\Bigg)\Bigg( F^+\left(u_1-\frac{c}{4}\right)-F^-\left(u_1+\frac{c}{4}\right)\Bigg)
\en
to yield the required result. 

(iv).
Finally,  recalling the identity given in \eqref{etadelta}, the $[E(u_1),F(u_2)]$ relation is proven by substituting  \eqref{rll1}[[2,4]] in \eqref{rll1}[[2,3]] to obtain:
\bea
&&\left[E^\pm (u_1),F^\pm (u_2)\right] = q^{2u}\left(q^{-1}-q\right)\Big( K_2^\pm (u_2)^{-1}K_1^\pm (u_2)\frac{\eta (P-u-1)}{\eta (u)\eta (P-1)} \nn\\
&&\hspace{3.5cm}-K_2^\pm(u_1)^{-1}K_1^\pm(u_1)\frac{\eta (P+h-u-1)}{\eta (u)\eta (P+h-1)}\Big),\lb{plusminusrelef}
\ena
%\eqref{plusrelef}\footnote{ For the $[E^-(u_1),F^-(u_2)]$ term, replace $K_i^{+}\left(u\right)$ by $K_i^{-}\left(u\right)$ in \eqref{plusrelef}}
and by substituting  \eqref{rll2}[[2,4]] in \eqref{rll2}[[2,3]] to obtain:%\eqref{eq4.5}, \eqref{eq4.6} and \eqref{eq4.2}.
\bea  \lb{eplusfminuscomm}
&&\left[E^+(u_1),F^-(u_2)\right]=K_2^-(u_2)^{-1}K_1^-(u_2)\frac{q^{2(u-\frac{c}{2})-1}\eta(1)\eta (P-u+\frac{c}{2}-1)}{\eta (u-\frac{c}{2})\eta (P-1)} \nn\\
&& \hspace{1.5cm}-K_2^+(u_1)^{-1}K_1^+(u_1)\frac{q^{2(u+\frac{c}{2})-1}\eta(1)\eta (P-\frac{c}{2}+h-u-1)}{\eta (u+\frac{c}{2})\eta (P+h-1)},
\ena
\bea \lb{eminusfpluscomm}
&&\left[E^-(u_1),F^+(u_2)\right]=K_2^+(u_2)^{-1}K_1^+(u_2)\frac{q^{2(u+\frac{c}{2})-1}\eta(1)\eta (P-u-\frac{c}{2}-1)}{\eta (u+\frac{c}{2})\eta (P-1)}\nn\\
&& \hspace{1.5cm} -K_2^-(u_1)^{-1}K_1^-(u_1)\frac{ q^{2(u-\frac{c}{2})-1}\eta(1)\eta (P+\frac{c}{2}+h-u-1)}{\eta (u-\frac{c}{2})\eta (P+h-1)} ,
\ena
remembering to expand \eqref{eplusfminuscomm} about $z=0$ and \eqref{eminusfpluscomm} about $z=\infty$
\noi
This completes the verification that the map $\Phi$ in the theorem is an associative algebra homomorphism.  The $H$-algebra structure is defined the same way for both algebras.  It is easy to see that the $H$-bigrading and moment maps are preserved by $\Phi$, using the Gauss decomposition of $L^\pm(u)$ given in \eqref{gaussL} with Lemma \ref{halgUR} and the table following Lemma \ref {halgUR}.  Hence $\Phi$ is a $H$-algebra homomorhpism.  It is clearly surjective

The proof of the injectivity of $\Phi$ is analogous to \cite{DF} for $\phi : U_{q}(\glth) \longrightarrow U(R_q)$.  We view $U_q(\glth)$ and $U(R_q(\glth))$ as the degenerations as $x \to 0$ of $U_{q,x}(\glth)$ and $U(R)$, respectively.   Consider the  representation $(\pi_{\la,c}^x,V_{\la,c})$ of  $U_{q,x}(\glth)$ induced from an integrable representation $(\pi_{\la,c}^q,V_{\la,c})$ of  $U_{q}(\glth)$, with level $c$ (we can identify ${U_q(\glth)}\otimes \C[\hbs]$ with $U_q(\glth)$ and $U(R_q)\otimes \C[\hbs]$ with $U(R_q)$, because $P$-independence implies that $Q$ is central).  Since $\Phi$ is a homomorphism of $H$-algebras, we get the commutative diagram:
\\ 

\vspace{0.2cm}

\begin{tikzpicture}
  \matrix (m) [matrix of math nodes,row sep=4em,column sep=6em,minimum width=2em] {
     U_{q,x}/xU_{q,x} & U(R)/xU(R) \\
        {U_q}\otimes \C[\hbs] &U(R_q)\otimes \C[\hbs] \\
       {U_q} &U(R_q) & End(V_{\lambda,c})\\
     };
  \path[-stealth]

    (m-1-1) edge node [left] {$\simeq$} (m-2-1)
                edge node [above] {$\Phi$} (m-1-2)
         
   (m-2-1.east|-m-2-2) edge node [below] {$\simeq$} node [above] {$\phi \otimes Id$} (m-2-2)

     (m-2-1) edge node [left] {$\simeq$} (m-3-1)
     (m-1-2) edge node [left] {$\simeq$} (m-2-2)
      edge [dashed,->] node [right] {$\pi_{\lambda,c}^x$} (m-3-3) 
      % edge [double] node [right] {$\pi_x$} (m-3-2)
   (m-2-2) edge node [left] {$\simeq$}(m-3-2)
         
    (m-3-1.east|-m-3-2) edge node [below] {$\simeq$} node [above] {$\phi$} (m-3-2)
            %edge [dashed,->] node [right] {$\pi_{\lambda,c}^x$} (m-3-3) 
           % edge [double] node [right] {$\pi_x$} (m-3-2)
    (m-3-2) edge node [above] {$\pi_{\lambda,c}^q$}(m-3-3);
\end{tikzpicture}

\noi
%%%IS THIS NECESSARY??
  Hence $\ker \Phi \subseteq x U_{q,x}(\glth)$.  Clearly, $\ker \Phi \subseteq \bigcap_{\la} \ker \pi_{\la,c}^x \circ \Phi$. Now, it is well-known from the representation theory of $U_q(\widehat{gl_n})$, that $\bigcap_{\la,c} \ker \pi_{\la,c}^q=\left\{0\right\}$ (see eq(5.73) in \cite{DF}).  Since $\ker \pi_{\la,c}^x \subseteq \ker \pi_{\la,c}^q$, it follows that  $\ker \Phi =\left\{0\right\}$, which completes the verification that $\Phi$ is injective.
\qed\\

\noi
{\em Remark.} The pair of relations \eqref{eplusfminuscomm} and  \eqref{eminusfpluscomm} can be directly verified (the proof is similar to \ref{plusrelef}). This indicates their consistency with the Definition \ref{halfcurrentseries} of (negative) half-currents. 
%%%%%%%%%%%%%%%%%%%%%%%%%%%%%%%%%%%%%%%%%%%%%%%%%%%%
\noi

\section{$H$-Hopf Algebroid Structures}
Let $A$ be a $H$-algebra (\ref{halg}).  We now recall the definition of $H$-bialgebroid and $H$-Hopf algebroid structures on $A$.  Start with the dynamical tensor product:
\begin{dfn}[Tensor Product of $H$-algebras $A\tot B$]
The tensor product of $A$ and $B$ denoted by $A \tot B$ is the $H^*$ bigraded vector space with
\be
 (A {\widetilde{\otimes}}B)_{\al\beta}=\bigoplus_{\gamma\in\h^*} (A_{\al\gamma}\otimes_{M_{\h^*}}B_{\gamma\beta}),
\en
where $\otimes_{M_{\hbs}}$ is the usual tensor product $\otimes$ modulo the relation:
\bea
&&\mu_r(f)a \tilde{\otimes } b = a \tilde{\otimes } \mu_l(f)b.
\ena
\end{dfn}
\noi
We will need a certain $H$-algebra of automorphisms that plays the role of the unit object:  \begin{dfn} [Algebra of Shift Operators $\mathcal{D}$]
Let $\mathcal{D}=\left\{  \sum _k {\hat{f}_kT_{\beta _k}\ } |\  \hat{f}_k \in M_{H^*}\ , \ \beta_k \in H^*  \right\}$. 
\end{dfn}
\noi
The bigrading and moment maps on $\cD$ are given as:
\be
&&\cD_{\al\al}=\left\{ \sum {\hat{f}T_{\alpha}}\ | \  \hat{f} \in M_{H^*}, \ \alpha \in H^*\right\}, \qquad \cD_{\al\beta}=0 \text{ if } \alpha \neq \beta\\
&&\mu_l(\hf)=\mu_r(\hf)=\hf T_{00}.
\en

\noi
{\em Remark.}  Note that $a\cong a\tot T_{-\beta}\cong T_{-\al}\tot a$ for all $a\in A_{\al\beta}$, proving the canonical isomorphism
\bea \lb{Dalg}
&&\mathcal{D} \tilde{\otimes} A\simeq \mathcal{D} \simeq  A\tilde{\otimes}\mathcal{D},
\ena
where the tensor product $\tpt$ is the usual $\otimes$ modulo the relation 
\bea \lb{dtprel}
&&f(u,P)a \tilde{\otimes } b = a \tilde{\otimes } f(u,P+h)b.
\ena

We will write the $R$ matrix elements in ${M}_{\h^*}$ as $(\widehat{R}_{(u)}^+)_{ab}^{xy}\equiv R^{+}(u,P)_{ab}^{xy}$ , where $a, b, x, y \in \left\{1,2\right\}$ (the inner subscript ${(u)}$ will be omitted when the context is clear).
We then have 
\bea \lb{rmatrixcoefficients}
&&\mu_l((\widehat{R}^+)_{ab}^{xy})=R^+(u,P+h)_{ab}^{xy},\quad  
\mu_r((\widehat{R}^+)_{ab}^{xy})
=R^{+}(u,P)_{ab}^{xy}.
\ena
The following definitions are well-known (see for example \cite{K2}):
\noi
\begin{dfn}[$H$-Bialgebroid]
An $H$-Hopf Algebroid $(A,H,\Delta,\varepsilon)$ is an $H$-algebra with the comultiplication and counit maps:
\be
&&\Delta:A \rightarrow A \tot A , \qquad \varepsilon: A \rightarrow \D,
\en
which are required to satisfy the compatibility conditions:
\bea \lb {ccepsilondelta}
&& (\Delta \tpt \id)\circ \Delta=(\id \tpt \Delta)\circ \Delta, \nn\\
&& (\varepsilon \tpt \id)\circ \Delta = (\id \tpt \varepsilon )\circ \Delta=\id.
\ena
\end{dfn}

\begin{dfn}[$H$-Hopf algebroid]
An $H$-bialgebroid $A$ is an $H$-Hopf algebroid $(A,H,\Delta,\varepsilon, S)$ if there is an algebra antihomomorphism $S:A \rightarrow A$ satisfying the compatibility conditions:
\bea \lb{ccantipode}
&& S(\mu_r(\hf))=\mu_l(\hf),\quad S(\mu_l(\hf))=\mu_r(\hf),\nn\\
&&m\circ (\id \otimes S)\circ\Delta(x)=\mu_l(\vep(x)1),\quad \forall x\in A,\\
&&m\circ (S\otimes\id  )\circ\Delta(x)=\mu_r(T_{\alpha}(\vep(x)1)),\quad \forall x\in A_{\alpha \beta}.
\ena
\end{dfn}

\subsection{$H$-Hopf Algebroid Structure on $U(R)$} \lb{hopfalgebroidUR} 
Let us define the $H$-bialgebroid maps for the $H$-algebra $U(R)$.
\begin{lem} \lb {bialg}
The Comultiplication $\Delta:U(R^\pm) \rightarrow U(R^\pm) \tot U(R^\pm)$ is given by  
\bea
&&\Delta \left(\mu _l(\hat{f})\right) = \mu _l(\hat{f})\tilde{\otimes }1,\qquad
\Delta \left(\mu _r(\hat{f})\right) = 1\tilde{\otimes }\mu _r(\hat{f}),\\
&&\Delta(e^Q)=e^Q \tot e^Q, \qquad \Delta \left(L^\pm_{\text{ab}}(u)\right)=\sum _{k=1}^2 L_{\text{ak}}^\pm(u)\tilde{\otimes }L_{\text{kb}}^\pm(u),
\ena
and the Counit $\varepsilon: U(R^\pm) \rightarrow \D$ is defined by
\bea 
 &&\varepsilon \left(\mu _l(\hat{f})\right)=\varepsilon \left(\mu _r(\hat{f})\right)=\hat{f}T_0, \lb{epsilonUR1}\\
&&\varepsilon \left(e^Q\right)=e^Q, \qquad \varepsilon \left(L_{\text{ab}}^\pm(u)\right)=\delta _{a,b}T_{\text{w(b)Q}\text{  }}. \lb{epsilonUR2}
\ena
\end{lem}
\noi
{\em {Proof}}.  

The verification of the compatibility conditions \eqref{ccepsilondelta} is straightforward, using the isomorphism \eqref{Dalg}.  We show the $\Delta$-invariance of the relations \eqref{rll1} (a similar calculation shows the invariance of the relations \eqref{rll2}, by using $c=c \tot 1 + 1 \tot c$).  
Let $R(u,P)(e_a \tot e_b) = \sum_{x,y}R(u,P)_{xy}^{ab}e_x \tot e_y$ and define the weight function $w:\left\{1,2\right\} \rightarrow \left\{\pm1\right\}$ by the identification:
\be
&&R(u,P)_{xy}^{ab} \rightarrow \left(
\begin{array}{cccc}
 R_{++}^{++} & R_{++}^{+-} & R_{++}^{-+} & R_{++}^{--} \\
 R_{+-}^{++} & R_{+-}^{+-} & R_{+-}^{-+} & R_{+-}^{--} \\
 R_{-+}^{++} & R_{-+}^{+-} & R_{-+}^{-+} & R_{-+}^{--} \\
 R_{--}^{++} & R_{--}^{+-} & R_{--}^{-+} & R_{--}^{--} \\
\end{array}
\right).\nn\\
\en
 Then the relations in \eqref{rll1} become:
\be
\sum_{a,b}R^\pm(u,P+h)_{cd}^{ab} \ L^\pm_{aa'}(u_1)\ L^\pm_{bb'}(u_2)&=&\sum_{c', d'}
L^\pm_{dc'}(u_2) \ L^\pm_{cd'}(u_1) \ R^{\pm}(u,P)_{d'c'}^{a'b'}.
\en
Apply $\Delta$ on both sides.
\be
\Delta(LHS)&=&\sum_{a,b}
\Delta(R^\pm(u,P+h)_{cd}^{ab})
\Delta(L^\pm_{aa'}(u_1))\Delta(L^\pm_{bb'}(u_2))\\
&=&\sum_{a,b \atop c', d'}
R^\pm(u,P+h)_{cd}^{ab}L^\pm_{ad'}(u_1)L^\pm_{bc'}(u_2)\widetilde{\otimes} L^\pm_{d'a'}(u_1)L^\pm_{c'b'}(u_2)\\
&=&\sum_{a',b'\atop c', d'}
L^\pm_{db}(u_2)L^\pm_{ca}(u_1)
R^{\pm}(u,P)_{ab}^{d'c'}\widetilde{\otimes} L^\pm_{d'a'}(u_1)L^\pm_{c'b'}(u_2),
\en%&=&\sum_{a',b'\atop c', d'}R^{+*}(u,P)_{ab}^{d'c'}
 \be
\Delta(RHS)&=&\sum_{c', d'}
\Delta(L^\pm_{dc'}(u_2))\Delta(L^\pm_{cd'}(u_1))\Delta(R^{\pm}(u,P)_{d'c'}^{a'b'} )\\
&=&\sum_{a,b\atop c', d'}
L^\pm_{db}(u_2)L^\pm_{ca}(u_1)
\widetilde{\otimes} L^\pm_{bc'}(u_2)L^\pm_{ad'}(u_1)
R^{\pm}(u,P)_{d'c'}^{a'b'}\\
&=&\sum_{a',b'\atop c', d'}
L^\pm_{db}(u_2)L^\pm_{ca}(u_1)
\widetilde{\otimes} R^{\pm}(u,P+h)_{a'b'}^{d'a'}L^\pm_{d'a'}(u_1)L^\pm_{c'b'}(u_2)\\
&=&\sum_{a',b'\atop c', d'}R^{\pm}(u,P)_{ab}^{d'c'}
L^\pm_{db}(u_2)L^\pm_{ca}(u_1)
\widetilde{\otimes} L^\pm_{d'a'}(u_1)L^\pm_{c'b'}(u_2).
\en
The final equality was obtained by using the following formula:
\bea
&&f(u,P)a\tot b=a\tot f(u,P+h)b\qquad a,b\in U_{q,x}(\glth),\lb{fstotf}
\ena
The equality of both expressions $\Delta(LHS)$ and $\Delta(RHS)$ follows from the fact:
\be
&&R^{\pm}(u,P+w(a)+w(b))_{ab}^{d'c'}=R^{\pm}(u,P)_{ab}^{d'c'}.
\en
\noi
The invariance of the first pair of  $RLL$-relations \eqref{rll1} under $\varepsilon$ follows using \eqref{epsilonUR1},  \eqref{epsilonUR2} and \eqref{rmatrixcoefficients}. 
\qed
%%antipode

Let us define the antipodal map using the inverse of $L(u)$.  It coincides with the elliptic case at $p=0$. 
\begin{lem} \lb{antipode}
 For $L^\pm(u)=((L^\pm_{ij}(u)))_{i,j=1,2}$, the antipode $S: U(R)\rightarrow U(R)$ is an algebra antihomorphism given by
\bea 
&&S(e^Q)=e^{-Q},\qquad S(\mu_r(\hf))=\mu_l(\hf),\quad S(\mu_l(\hf))=\mu_r(\hf),\nn\\
&&S(L^\pm_{11}(u))=L^\pm_{22}(u-1), \qquad \qquad \qquad S(L^\pm_{12}(u))=-q^{-1}\frac{\eta(P+h+1)}{\eta(P+h)}
L^\pm_{12}(u-1), \quad \nn\\
&&S(L^\pm_{21}(u))=-q\frac{\eta(P)}{\eta(P+1)}L^\pm_{21}(u-1),\quad S(L^\pm_{22}(u))=\frac{\eta(P+h+1)\eta(P)}{\eta(P+h)\eta(P+1)}\L^\pm_{11}(u-1).\nn
\ena
\end{lem}
\noi
{\em {Proof}}.   For $i, j = 1, 2$, by definition  $S(L_{ij}^\pm(u))=(L^\pm(u)^{-1})_{ij}$.  Then, the fact that the $RLL$- relations are invertible implies their $S$-invariance.  For the calculation that the map $S$ is compatible with $\varepsilon$ and $\Delta$,  we use the expansion of the $RLL$-relations in the appendix.   Let us verify \eqref{ccantipode} for $L^+_{12}(u)$.  The relation 
\eqref{alphabeta} with $u_1=u-1, u_2=u$ and $u=-1$ becomes
\bea
L_{11}^{\pm}(u)L_{12}^{\pm}(u-1) = \frac{q\eta(P-1)}{\eta(P)}L_{12}^\pm(u)L_{11}^\pm(u-1),
\ena
which yields the final equality in the following calculation:

\be %epsdelSproof
&&m\circ (id\otimes S  )\circ\Delta(L^\pm_{12}(u))\nn\\
&&\quad =-L_{11}^{\pm}(u) \frac{q^{-1}\eta(P+h+1)}{\eta(P+h)}L_{12}^{\pm}(u-1)+L_{12}^\pm(u)\frac{\eta(P+h+1)\eta(P)}{\eta(P+h)\eta(P+1)}L_{11}^\pm(u-1)\nn\\
&&\quad  = \frac{q^{-1}\eta(P+h+2)}{\eta(P+h+1)}L_{11}^{\pm}(u)L_{12}^{\pm}(u-1) + \frac{\eta(P+h+2)\eta(P-1)}{\eta(P+h+1)\eta(P)}L_{12}^\pm(u)L_{11}^\pm(u-1)\nn\\
&&\quad =0.
\en
Since $\varepsilon(L^+_{12}(u)=0$, the relation \eqref{ccantipode} is confirmed.  The proofs for the remaining generators $L^\pm_{11}(u), L^\pm_{21}(u)$ and $L^\pm_{22}(u)$ are similar.
\qed

\begin{thm} \lb{thmUR}
The following statements are true:
\begin{enumerate}
\item The $H$-algebras $U(R^{\pm}):= U\left(R^\pm_{q,x}(\glth)\right)$ are $H$-Hopf Algebroids.
\item The Total $H$-algebra $U(R)$ is an $H$-bialgebroid.
\end{enumerate}
\end{thm}
\noi
{\em {Proof.}} 
(i). The $H$-bialgebra structure is proven in Lemma \ref{bialg}, while the antipodal map is established  in Lemma \ref{antipode}.
Statement (ii) follows from (i) and one easily checks that the second pair is $\Delta$-invariant using $\Delta(q^c)=q^c \tot q^c$.
\qed
%%%
%%Determinant
%%%

In order to develop the representation theory and to make contact with the more familiar dynamical algebras, we will consider the subalgebras $U_{q,x}(\slth)$ and $U(R_{q,x}(\slth))$ of  $U_{q,x}(\glth)$ and $U(R_{q,x}(\glth))$ respectively.
The $H$-Hopf algebroid structures is also more transparent on these subalgebras.  We begin by introducing the dynamical determinant element, similar to the finite-dimensional case \cite {KR} as follows:

\begin{dfn} \lb{det}
The Dynamical Determinant element is defined as
\be 
&& Det(L^{\pm}(u))=L^{\pm}_{11}(u+1)L^{\pm}_{22}(u) - \frac{q\eta(P-1)}{\eta(P)} L^{\pm}_{12}(u+1)L^{\pm}_{21}(u).
\en
\end{dfn}
\noi
The main properties are given below:
\begin{prop}
The element $Det(L^{\pm}(u))$ is central in $U(R)$.  Further,
\bea \lb{deltadet}
&&\Delta(Det(L^{\pm}(u)))=Det(L^{\pm}(u)) \tot Det(L^{\pm}(u)), \quad \varepsilon(Det(L^{\pm}(u)))=1.
\ena
\end{prop}
% dr k2
%%

\noi
{\em Proof}.  Using Theorem \ref{mainthm}, we can apply the Gauss decomposition  \eqref{gaussL} of $L^{\pm}(u)$ along with the commutation relation \eqref{plusrelek2} to obtain 
\bea
&&Det(L^{\pm}(u))=K^{\pm}_1(u+1)K^{\pm}_2(u).
\ena  Then it is straightforward to verify that this element lies in the center of $U(R)$ (resp. $U(R^+))$ by using the defining relations in Proposition \ref{Defrelns} (resp. \ref{plusrelns}) and the next expression which is obtained by expanding \eqref{rho}
\bea \lb{rholid}
&&\rho^+(u)\rho^+(u+1)=q^{-1}\hspace{0.1cm}\frac{\eta(u+1)}{\eta(u)}.  
\ena
For the first expression in \eqref{deltadet}, expand both the sides using Definition \ref{det} and the coproduct in Lemma \ref{bialg}.  Now use the formulae
\bea
&&\frac{\bar{c}(1,P)}{\bar{b}(1,P)}=\frac{b(1,P)}{c(1,P)}=q \frac{\eta(P-1)}{\eta(P)}
\ena
in relations \eqref{alphabeta} and \eqref{A5}, with $u_1=u+1$ and $u_2=u$.   It remains to show that
\bea \lb{detproof}
&&L^{\pm}_{12}(u+1)L^{\pm}_{21}(u) \tot D(u)=0
\ena
where $D(u)$ is given by
\bea
&&D(u)=q\frac{\eta(P+h-1)}{\eta(P+h)}L^{\pm}_{11}(u+1)L^{\pm}_{22}(u)-q\frac{\eta(P-1)}{\eta(P)}L^{\pm}_{22}(u+1)L^{\pm}_{11}(u) \nn\\
&&\qquad \qquad + q\frac{\eta(P-1)\eta(P+h-1)}{\eta(P)\eta(P+h)}L^{\pm}_{12}(u+1)L^{\pm}_{21}(u)-L^{\pm}_{11}(u+1)L^{\pm}_{22}(u).
\ena
\eqref{detproof} can be checked by multiplying \eqref{A.11} and \eqref{A.12} (with $u_1=u+1$ and $u_2=u$) by $b(1,P)\bar{c}(1,P)$ and $c(1,P)\bar{c}(1,P)$ respectively, and subtracting them.  The second formula in \eqref{deltadet} is an easy consequence of the counit definition in Lemma \ref{bialg}.
\qed 

\begin{dfn}
The subalgebra $U(R_{q,x}(\slth))$ of $U(R_{q,x}(\glth))$ is defined by the condition:
\be
&&Det(L^{\pm}(z))=1. 
\en
\end{dfn}
%%%%
%%%%%%
%%%%%%%%%%%%%%
\noi
We will obtain a similar statement as Theorem \ref{mainthm} for $\slth$
in the next subsection.
\subsection {Relation to Standard Drinfeld Realization of \texorpdfstring{$U_{q}(\slthBig)$}{quantum affine sl2}}\lb{drinfeld}
%%%%
Consider a field $\K \supseteq \C$. The following standard presentation, known as the Drinfeld realization, is well-known \cite{Drinfeld}.

\begin{dfn}\lb{defstandard}
$\K[U_q(\slth)]$ is the associative algebra over $\K$ 
generated by 
$a_n\ (n\in \Z_{\not=0})$, $x_n^\pm\ (n\in \Z)$, 
$h$,  $c$ and $d$. 
The defining relations are given as follows.
\be
&&c :\hbox{ central },\nn\\
&& [h,d]=0,\quad [d,a_{n}]=n a_{n},\quad 
[d,x^{\pm}_{n}]=n x^{\pm}_{n}, \nn\\
&&[h,a_{n}]=0,\qquad [h, x_n^\pm]=\pm 2 x_n^\pm,\nn\\
&&
[a_{m},a_{n}]=\frac{[2n]_{q}[c n]_{q}}{n}\delta_{n+m,0},\nn
\\&&
[a_{n},x_m^+]=\frac{[2n]_{q}}{n}q^{\frac{c|n|}{2}}x_{m+n}^+,\nn\\
&&[a_{n},x_m^-]=-\frac{[2n]_{q}}{n}q^{-\frac{c|n|}{2}} x_{m+n}^-,\nn\\
&&x_{m+1}^\pm x_n^\pm - q^{\pm 2}x_n^\pm x_{m+1}^\pm =  q^{\pm 2}x_m^\pm x_{n+1}^\pm - x_{n+1}^\pm x_m^\pm ,\nn\\
%endo
&&[x^+_n,x^-_m] =\frac{1}{q-q^{-1}}\left(q^{\frac{c(m-n)}{2}}\psi _{m+n}-q^{\frac{c(n-m)}{2}}\varphi _{m+n}\right). \nn
\en

Denote $x^\pm(z)=\sum_{n\in \Z}x^\pm_{n} z^{-n}$, and
the auxillary currents $\psi_n$,\ $\varphi_{-n}$, ($n \geq 0$) by
\be
&&\sum_{n\geq 0}\psi_n z^{-n}=q^{h}
\exp\left( (q-q^{-1}) \sum_{n>0} a_{n}z^{- n}\right),\quad 
\sum_{n\geq 0}\varphi_{-n}z^n=q^{-h}
\exp\left(-(q-q^{-1})\sum_{n>0} a_{-n}z^{ n}\right).\\
\en

\end{dfn}

Choosing $\K=\C[H]$ and defining $\widehat U_{q,x}(\slth):=\K[U_q(\slth)] \otimes \C[\hbs]$,
consider the two operators {\footnote{$K^\pm(u)$ can be obtained as suitable degenerations of the expressions for the elliptic operators $K(z)$ and $H^\pm(u)$ in eqns (3.15, 3.25, 3.29) \cite{JKOS2} (without the elliptic shift by the central element: $p \to pq^{-2c}$).}}
\bea \lb{defnKplus}
&& K^+(u)=\exp\left(\sum_{n > 0} \frac{[n]_q}{[2n]_q}(q-q^{-1})a_n q^{-(2u+1)n}\right)  e^Q q^{\frac{h}{2}},\\
&& K^-(u)=\exp\left(-\sum_{n > 0} \frac{[n]_q}{[2n]_q}(q-q^{-1})a_{-n} q^{(2u+1)n}\right)  e^Q q^{-\frac{h}{2}}.
\ena
The $H$-algebra structure on  $\widehat U_{q,x}(\slth)$ is defined in exactly the same way as  $U_{q,x}(\glth)$, using the same Heisenberg algebra $\H$ in Subsection \ref{heis}.
Now we define the dynamical currents by
\bea
&& \widehat{K_1^\pm}(u)=K^{\pm}(u-1), \quad \widehat{K_2^\pm}(u)=K^{\pm}(u)^{-1}, \\
&&\widehat{E}(u)=x^+(z)e^{2Q}, \quad \widehat{F}(u)=x^-(z), \quad \widehat{H^{\pm }}(u)=\widehat{K_1^{\pm }}(u)\widehat{K_2^{\pm }}(u)^{-1}. \lb{hatefhdefn}
\ena
Define the derived subalgebra  $\widehat U_{q,x}(\slth)'$ of $\widehat U_{q,x}(\slth)$ by the same relations as Definition \ref{defstandard}, but without the element $d$.  We get the following  {Corollary to Theorem \ref {mainthm}}:
\begin{cor}\lb {mainthmsl2} The isomorphism in Theorem \ref{mainthm} restricts to an $H$-subalgebra isomorphism:
$\widehat{U}_{q,x}(\slth)' \simeq  U(R_{q,x}(\slth)).$
\end{cor}
\noi
{\em Proof.}  We identify $\widehat{U}_{q,x}(\slth)'$ with the $H$-algebra generated by
\be
&& \left\{ h,\ c,\ \widehat{E}(u),\ \widehat{F}(u),\ \widehat{K_1^\pm}(u),\  \widehat{K_2^\pm}(u),\  \widehat{K_1^\pm}(u)^{-1},\  \widehat{K_2^\pm}(u)^{-1}\right\}.
\en

\noi
Replacing $\widehat{E(u)}$ with $(q-q^{-1})\widehat{E(u)}$ and  $\widehat{F(u)}$ with $(q-q^{-1})\widehat{F(u)}$, we can easily verify the defining relations of the algebra $U_{q,x}(\glth)$ in Proposition \ref{Defrelns}.  The calculation is using the Baker-Campbell-Hausdorff formula.  The relations between the $\widehat{K^\pm_i}(u)$ are derived from the relation for $[a_n,a_m]$ and Definition \ref{defnKplus} using the formula:
$ [A,B]=k \ \implies \ \exp(A)\exp(B)\exp(-A)=e^k \exp(B)$.

The relations between $\widehat{K^\pm_i}(u)$ and  $\left\{\widehat{E}(u), \widehat{F}(u) \right\}$ are verified by expressing the defining relations $\left[a_n,x^\pm_m\right]$ as
\be
&& [a_{n},x^+(z)]=\frac{[2n]_{q}}{n}q^{\frac{c|n|}{2}}z^n x^+(z),\qquad
 [a_{n},x^-(z)]=-\frac{[2n]_{q}}{n} q^{-\frac{c|n|}{2}}z^n x^-(z).
\en

\noi
and applying the fact:
\be
&& [A,X]=kX \implies \exp(A) X \exp(-A)=e^kX.
\en
It follows that  $q^hx^\pm(z)q^{-h}=q^{\pm 2}x^{\pm}(z)$.
\noi
Finally, the relation for $[\widehat{E}(u), \widehat{F}(u)]$ is a consequence of the defining relations in \eqref{hatefhdefn}, since   $[(q-q^{-1})x_n^+,(q-q^{-1})x_m^-]$ identifies with \eqref{efcomm}.
Thus  $\widehat{U}_{q,x}(\slth)'$ is a subalgebra of  $U_{q,x}(\glth)$.  The $L$-operators constructed by replacing $\left\{E^\pm(u), F^\pm(u), K_i^{\pm}(u) \right\}$ by $\left\{ \widehat{E}^\pm(u), \widehat{F}^\pm(u), \widehat{K^\pm_i}(u) \right\}$ are in $U(R^\pm_{q,x}(\slth))$ because 
\bea
&& Det(L^\pm(u))=\widehat {K^{\pm}_1}(u+1)\widehat{K^{\pm}_2}(u)=K^{\pm}(u)K^{\pm}(u)^{-1}=1,
\ena
where the dynamical determinant $Det$ is defined in Definition \ref{det}.
\qed

\noi
We will hereafter write  $U_{q,x}(\slth)$ for $\widehat{U}_{q,x}(\slth)'$ and the half-current subalgebras are given by
\be
U^{\pm}_{q,x}(\slth):=U_{q,x}(\slth) \cap U^{\pm}_{q,x}(\glth) \subseteq U_{q,x}(\slth).
\en
\subsection{$H$-Hopf algebroid structure on  \texorpdfstring{$U_{q,x}(\slthBig)$}{Uqx}}

%The {$H$-Hopf algebroid structure on  $U_{q,x}(\glth)$ can be restricted to $\slth$.
\begin{thm} \lb{thmuq}
The following statements are true.
\begin{enumerate}
\item
The Half-current Algebras $(U^\pm_{q,x}(\slth),\hbs, \mu_l, \mu_r, \Delta,\varepsilon,S)$ are $H$-Hopf-algebroids.
% called the Negative Dynamical Affine Trignometric Quantum Group,   NAQG.\item
\item
The Total Algebra $(U_{q,x}(\slth),\hbs,\mu_l,\mu_r)$ is an $H$-bialgebroid.
%called the Dynamical Affine Trignometric Quantum Algebra, DATQA.  Or we prefer Dynamical Quantum Universal Enveloping Algebra DQUEA.
\end{enumerate}
%The statements (i) and (ii) are also true for $U(R)$ and $U(R^\pm)$ respectively.
\end{thm}

\noi
{\em {Proof.}} (i) Since $U(R^\pm)$ is an $H$-Hopf algebroid (Subsection \ref{hopfalgebroidUR}), we can define the comultiplication, counit and antipodal maps on $U^\pm_{q,x}(\slth)$, by using Theorem \ref{mainthm} and the Gauss decomposition of  $L^\pm(u)$ given in \eqref{gaussL}.  The coproduct is given explicitly in Proposition \ref{comult}.  Use the expression \eqref{linv} for the inverse of $L^+(u)$ along with the commutation relations \eqref{plusrelek2} and \eqref{plusrelfk2} (with $u_1=u-1$), in $S(L^+(u))=L^+(u)^{-1}$ to obtain the explicit formula for the image of the half-currents under the antipodal map.  (ii) The statement follows from Theorem \ref{thmUR}.
\qed
 
Explicit expressions for the comultiplication map of the positive and negative half-currents are available in our case.  The elliptic version of the next result appears in Proposition 3.12 of \cite{K2}.
\begin{prop}\lb{comult}
The coproduct for the half-currents is given as:
\be 
&&\Delta(K_1^{\pm}(u))=K_1^{\pm}(u)\tot K_1^{\pm}(u)+\sum_{j=1}^\infty (-1)^jK_1^{\pm}(u)E^{\pm}(u-1)^j\tot 
F^{\pm}(u-1)^jK_1^{\pm}(u),\\
&&\Delta(K_2^{\pm}(u))=K_2^{\pm}(u) \tot K_2^{\pm}(u) + K_2^{\pm}(u)  E^{\pm}(u) \tot F^{\pm}(u) K_2^{\pm}(u),\\
&&\Delta(E^{\pm}(u))=1\tot E^{\pm}(u)+E^{\pm}(u)\tot K_2^{\pm}(u)^{-1}K_1^{\pm}(u)\\
&&\qquad\qquad\qquad\qquad+\sum_{j=1}^\infty(-1)^jE^{\pm}(u)^{j+1}\tot K_2^{\pm}(u)^{-1}F^{\pm}(u)^jK_1^{\pm}(u),\\
&&\Delta(F^{\pm}(u))=F^{\pm}(u)\tot 1+ K_1^{\pm}(u)K_2^{\pm}(u)^{-1}\tot F^{\pm}(u)\\
&&\qquad\qquad\qquad\qquad+\sum_{j=1}^\infty(-1)^jK_1^{\pm}(u)E^{\pm}(u)^jK_2^{\pm}(u)^{-1}\tot F^{\pm}(u)^{j+1},\\
&& \Delta \left(H^{+}(u)\right)=H^{+}(u)\hat{\otimes } H^{+}(u) \bmod A_{\geqslant 0}\tpt A_{\leqslant 0},
\en
where $A_{\geqslant 0}$  (resp. $A_{\leqslant0}$) is the subalgebra generated by
 $\left\{q^{\pm c}, K_i^+(u), \ E^+(u)\ (\text{resp. } F^+(u)) \right\}.$
\end{prop}
\noi{\it Proof.} Straightforward, by using the coproduct formula in Proposition \ref{bialg} along with the Gauss decomposition of $L^\pm(u)$ in \eqref{gaussL} for the first four formulae.  The last relation is obtained by substituting the first two in the definition of $H^+(u)$. \qed 

 Thus, as suggested in \cite{K2} (without a proof), the Hopf algebroid structure on the     {\em elliptic} quantum group $U_{q,p}(\slth)$ does indeed survive the degeneration as $p \to 0$. 
 
\subsection{Concluding Remarks and Questions}
% The infinite dimensional representations of $U_{q,p}{\widehat{(\frak{g})}}$ in terms of dynamical quantum $Z$- and $W$- algebras are explicitly constructed.   However, the half-current algebras are not investigated there.
The main theorem in this article has been extended by the author to  $U_{q,x}(\glnh)$ and to $U_{q,x}(gl_n)$ and will appear elsewhere.  After completing this study, there appeared a new definition of elliptic $U_{q,p}{\widehat{(\frak{g})}}$ associated to any untwisted affine Lie algebra \cite{Knew} for arbitrary values of $p$, that can be identified at $p=0$ with the definition given here of $U_{q,x}(\slth)$ (see Theorem 2.2 in \cite{Knew}).   Our main result also extends to the elliptic case, i.e. to $U_{q,p}(\glnh)$ and will be discussed in a separate publication. Finally, it would be interesting to find a relationship of $U_{q,x}(\slth)$, when c=0, with the elliptic algebra $U(2)$ in \cite{KNR} .
%\pagebreak
\subsection*{Acknowledgments}
%\vspace{-37cm}
The author would like to thank Pavel Etingof for introducing him to the dynamical quantum groups and Hitoshi Konno for some valuable comments.   He also expresses his gratitude to  Nick Early, David Eisenbud, David Saint John, and Suresh Srinivasamurthy for their kind encouragement and advice.  This work was completed at MSRI, Berkeley and at Pennsylvania State University.   
 
%%%%%%%%%%%%%%%%%%%%%%%%%%%%%%%%%%%%%%%%%%%%%%%%%%%%%%%%%%%%%%%%%%%%%%%%%%%%%%%%%%%%%

\appendix
\noi 

\setcounter{equation}{0}
\begin{appendix}
 %%%%%%%%%%%%%%%%%%%%%%%%%%%%%%%
\pagebreak

\newpage
\section {The $RLL$ relations of $U(R)$}
\begin{prop} \lb{rllexpanded} Using the notation in \eqref{rappendix}, the first relation  \eqref{rll1}
\be
&&R^{\pm(12)}(u,P+h)L^{\pm (1)}(u_1)L^{\pm(2)}(u_2)=L^{\pm (2)}(u_2)L^{\pm (1)}(u_1)R^{\pm(12)}(u,P),
\en
 is expanded as
%{(ii). The $R^+L^+L^-$ relations}
\bea
&&L^{\pm}_{11} (u_1) L^{\pm}_{11} (u_2)=L^{\pm}_{11} (u_2) L^{\pm}_{11} (u_1),\quad
L^{\pm}_{12} (u_1) L^{\pm}_{12} (u_2)=L^{\pm}_{12} (u_2) L^{\pm}_{12} (u_1),\\
&&L^{\pm}_{21} (u_1) L^{\pm}_{21} (u_2)=L^{\pm}_{21} (u_2) L^{\pm}_{21} (u_1) ,\quad
L^{\pm}_{22} (u_1) L^{\pm}_{22} (u_2)=L^{\pm}_{22} (u_2) L^{\pm}_{22} (u_1),\\
&&L^{\pm}_{11} (u_1) L^{\pm}_{12} (u_2)=L^{\pm}_{11} (u_2) L^{\pm}_{12} (u_1) \bar{c}(u,P)+L^{\pm}_{12} (u_2) L^{\pm}_{11} (u_1) b(u,P),\lb{alphabeta}\\ 
&&L^{\pm}_{12} (u_1) L^{\pm}_{11} (u_2)=L^{\pm}_{11} (u_2) L^{\pm}_{12} (u_1) \bar{b}(u,P)+L^{\pm}_{12} (u_2) L^{\pm}_{11} (u_1) c(u,P),\\
&&L^{\pm}_{21} (u_1) L^{\pm}_{22} (u_2)=L^{\pm}_{21} (u_2) L^{\pm}_{22} (u_1) \bar{c}(u,P)+L^{\pm}_{22} (u_2) L^{\pm}_{21} (u_1) b(u,P),\lb{A5}\\
&&L^{\pm}_{22} (u_1) L^{\pm}_{21} (u_2)=L^{\pm}_{21} (u_2) L^{\pm}_{22} (u_1) \bar{b}(u,P)+L^{\pm}_{22} (u_2) L^{\pm}_{21} (u_1) c(u,P), \\
&&b(u,P+h) L^{\pm}_{11} (u_1) L^{\pm}_{21} (u_2)+ c(u,P+h) L^{\pm}_{21} (u_1) L^{\pm}_{11} (u_2)=L^{\pm}_{21} (u_2) L^{\pm}_{11} (u_1), \nn\\ \\
&&b(u,P+h) L^{\pm}_{12} (u_1) L^{\pm}_{22} (u_2)+ c(u,P+h) L^{\pm}_{22} (u_1) L^{\pm}_{12} (u_2)=L^{\pm}_{22} (u_2) L^{\pm}_{12} (u_1), \nn\\ \\
&&\bar{b}(u,P+h) L^{\pm}_{21} (u_1) L^{\pm}_{11} (u_2)+\text{  }\bar{c}(u,P+h) L^{\pm}_{11} (u_1) L^{\pm}_{21} (u_2)=L^{\pm}_{11} (u_2) L^{\pm}_{21} (u_1), \nn\\ \\
&&\bar{b}(u,P+h) L^{\pm}_{22} (u_1) L^{\pm}_{12} (u_2)+\text{  }\bar{c}(u,P+h) L^{\pm}_{12} (u_1) L^{\pm}_{22} (u_2)=L^{\pm}_{12} (u_2) L^{\pm}_{22} (u_1), \nn\\ \\
&&b(u,P+h) L^{\pm}_{11} (u_1) L^{\pm}_{22} (u_2)+c(u,P+h) L^{\pm}_{21} (u_1) L^{\pm}_{12} (u_2) \nn\\
&&\qquad \qquad =L^{\pm}_{21} (u_2) L^{\pm}_{12} (u_1) \bar{c}(u,P)+L^{\pm}_{22} (u_2) L^{\pm}_{11} (u_1) b(u,P),\lb{A.11}\\
&&b(u,P+h) L^{\pm}_{12} (u_1) L^{\pm}_{21} (u_2)+c(u,P+h) L^{\pm}_{22} (u_1) L^{\pm}_{11} (u_2) \nn\\
&&\qquad \qquad =L^{\pm}_{21} (u_2) L^{\pm}_{12} (u_1) \bar{b}(u,P)+L^{\pm}_{22} (u_2) L^{\pm}_{11} (u_1) c(u,P),\lb{A.12}\\
&&\bar{b}(u,P+h) L^{\pm}_{21} (u_1) L^{\pm}_{12} (u_2)+\bar{c}(u,P+h) L^{\pm}_{11} (u_1) L^{\pm}_{22} (u_2) \nn\\
&&\qquad \qquad =L^{\pm}_{11} (u_2) L^{\pm}_{22} (u_1) \bar{c}(u,P)+L^{\pm}_{12} (u_2) L^{\pm}_{21} (u_1) b(u,P),\lb{A.13}\\
&&\bar{b}(u,P+h) L^{\pm}_{22} (u_1) L^{\pm}_{11} (u_2)+\bar{c}(u,P+h) L^{\pm}_{12} (u_1) L^{\pm}_{21} (u_2) \nn\\
&&\qquad \qquad =L^{\pm}_{11} (u_2) L^{\pm}_{22} (u_1) \bar{b}(u,P)+L^{\pm}_{12} (u_2) L^{\pm}_{21} (u_1) c(u,P). \lb{A.14} 
\ena
%{(iii). The inverse of the R-matrix}
The relation \eqref{rll2}
\be
&&R^{\pm(12)}\left(u\pm\frac{c}{2},P+h\right)L^{\pm(1)}(u_1)L^{\mp(2)}(u_2)=L^{\mp(2)}(u_2)L^{\pm(1)}(u_1)R^{\pm(12)}\left(u\mp\frac{c}{2},P\right),
\en
\vspace{-20pt}
is expanded below.  We write the $R$-matrix entries $c(u,P)$ and $\bar{c}(u,P)$ (see  \eqref{rappendix}) as $c_0(u,P)$ and $\bar{c_0}(u,P)$ respectively,  to distinguish them from the central element $c$.
\bea
&& \rho^{\pm}(u \pm \frac{c}{2}) \ L^{\pm}_{11}(u_1)  L^{\mp}_{11}(u_2)= L^{\mp}_{11}(u_2)  L^{\pm}_{11}(u_1) \ \rho^{\pm}(u \mp \frac{c}{2}) ,\\
&& \rho^{\pm}(u \pm \frac{c}{2}) L^{\pm}_{11}(u_1)  L^{\mp}_{12}(u_2) \nn\\
&&\qquad \quad  =\Big( L^{\mp}_{11}(u_2)  L^{\pm}_{12}(u_1) \bar{c_0}(u \mp \frac{c}{2},P)+ \  L^{\mp}_{12}(u_2) L^{\pm}_{11}(u_1) b(u \mp \frac{c}{2},P)\Big) \rho^{\pm}(u \mp \frac{c}{2}),\nn\\ \\
&& \rho^{\pm}(u \pm \frac{c}{2})  L^{\pm}_{12}(u_1)  L^{\mp}_{11}(u_2) \nn\\
&& \qquad \quad  = \Big(L^{\mp}_{11}(u_2)  L^{\pm}_{12}(u_1) \bar{b}(u \mp \frac{c}{2},P)+ \  L^{\mp}_{12}(u_2) L^{\pm}_{11}(u_1) c_0(u \mp \frac{c}{2},P) \Big)\rho^{\pm}(u \mp \frac{c}{2}),\nn\\ \\
&&\rho^{\pm}(u \pm \frac{c}{2})  L^{\pm}_{12}(u_1)  L^{\mp}_{12}(u_2)= L^{\mp}_{12}(u_2)  L^{\pm}_{12}(u_1) \rho^{\pm}(u \mp \frac{c}{2}),\\
&&\rho^{\pm}(u \pm \frac{c}{2})\Big(b(u \pm \frac{c}{2},P+h) L^{\pm}_{11}(u_1)  L^{\mp}_{21}(u_2) + c_0(u \pm \frac{c}{2},P+h)  L^{\pm}_{21}(u_1)  L^{\mp}_{11}(u_2)\Big)\nn\\  
&&\qquad \quad =L^{\mp}_{21}(u_2)  L^{\pm}_{11}(u_1)  \rho^{\pm}(u \mp \frac{c}{2}),\\
&&\rho^{\pm}(u \pm \frac{c}{2}) \Big( b(u \pm \frac{c}{2},P+h) L^{\pm}_{11}(u_1) L^{\mp}_{22}(u_2)+ c_0(u \pm \frac{c}{2},P+h)  L^{\pm}_{21}(u_1)  L^{\mp}_{12}(u_2)\Big)\nn\\
&&\qquad \quad  = \Big(L^{\mp}_{22}(u_2) L^{\pm}_{11}(u_1) b(u \mp \frac{c}{2},P)+ L^{\mp}_{21}(u_2)  L^{\pm}_{12}(u_1) \bar{c_0}(u \mp \frac{c}{2},P) \Big)\rho^{\pm}(u \mp \frac{c}{2}),\nn\\ \\
&&\rho^{\pm}(u \pm \frac{c}{2}) \Big( b(u \pm \frac{c}{2},P+h)   L^{\pm}_{12}(u_1)  L^{\mp}_{21}(u_2)+ c_0(u \pm \frac{c}{2},P+h) L^{\pm}_{22}(u_1)  L^{\mp}_{11}(u_2)\Big)\nn\\
&&\qquad \quad  = \Big(L^{\mp}_{21}(u_2)  L^{\pm}_{11}(u_1) {c_0}(u \mp \frac{c}{2},P) + L^{\mp}_{21}(u_2)  L^{\pm}_{12}(u_1) \bar{b}(u \mp \frac{c}{2},P)  \Big)\rho^{\pm}(u \mp \frac{c}{2}),\nn\\ \\
&&\rho^{\pm}(u \pm \frac{c}{2})\Big( b(u \pm \frac{c}{2},P+h)  L^{\pm}_{12}(u_1)  L^{\mp}_{22}(u_2)+c_0(u \pm \frac{c}{2},P+h)  L^{\pm}_{22}(u_1)  L^{\mp}_{12}(u_2)\Big)\nn\\
&&\qquad \quad = L^{\mp}_{22}(u_2)  L^{\pm}_{12}(u_1) \rho^{\pm}(u \mp \frac{c}{2}),\\
&&\rho^{\pm}(u \pm \frac{c}{2})\Big( \bar{b}(u \pm \frac{c}{2},P+h)  L^{\pm}_{21}(u_1)  L^{\mp}_{11}(u_2)+\bar{c_0}(u \pm \frac{c}{2},P+h) L^{\pm}_{11}(u_1)  L^{\mp}_{21}(u_2)\Big)\nn\\
&&\qquad \quad = L^{\mp}_{11}(u_2)  L^{\pm}_{21}(u_1) \rho^{\pm}(u \mp \frac{c}{2}),\nn\\
%\ena
%\bea
&& \rho^{\pm}(u \pm \frac{c}{2})\Big( \bar{b}(u \pm \frac{c}{2},P+h) L^{\pm}_{21}(u_1)  L^{\mp}_{12}(u_2)+\bar{c_0}(u \pm \frac{c}{2},P+h) L^{\pm}_{11}(u_1)  L^{\mp}_{22}(u_2)\Big)\nn\\
&&\qquad \quad  = \Big( L^{\mp}_{11}(u_2)  L^{\pm}_{22}(u_1) \bar{c_0}(u \mp \frac{c}{2},P) + L^{\mp}_{12}(u_2)  L^{\pm}_{21}(u_1) b(u \mp \frac{c}{2},P)  \Big)\ \rho^{\pm}(u \mp \frac{c}{2}),\nn\\ \\
&& \rho^{\pm}(u \pm \frac{c}{2})\Big(\bar{b}(u \pm \frac{c}{2},P+h)  L^{\pm}_{22}(u_1)  L^{\mp}_{11}(u_2)+\bar{c_0}(u \pm \frac{c}{2},P+h) L^{\pm}_{12}(u_1)  L^{\mp}_{21}(u_2) \Big) \nn\\
&&\qquad \quad=\Big(L^{\mp}_{11}(u_2)  L^{\pm}_{22}(u_1) \bar{b}(u \mp \frac{c}{2},P) + L^{\mp}_{12}(u_2)  L^{\pm}_{21}(u_1) {c_0}(u \mp \frac{c}{2},P) \Big) \rho^{\pm}(u \mp \frac{c}{2}),\nn\\ 
\ena
\bea
&& \rho^{\pm}(u \pm \frac{c}{2}) \Big(\bar{b}(u \pm \frac{c}{2},P+h) L^{\pm}_{22}(u_1)  L^{\mp}_{12}(u_2)+\bar{c_0}(u \pm \frac{c}{2},P+h) L^{\pm}_{12}(u_1) L^{\mp}_{22}(u_2)\Big)\nn\\
&&\qquad \quad = L^{\mp}_{12}(u_2)  L^{\pm}_{22}(u_1)\rho^{\pm}(u \mp \frac{c}{2}),\\
&&\rho^{\pm}(u \pm \frac{c}{2})  L^{\pm}_{21}(u_1)  L^{\mp}_{21}(u_2)= L^{\mp}_{21}(u_2)  L^{\pm}_{21}(u_1) \rho^{\pm}(u \mp \frac{c}{2}),\\
&&\rho^{\pm}(u \pm \frac{c}{2})  L^{\pm}_{21}(u_1)L^{\mp}_{22}(u_2)\nn\\
&&\qquad \quad = \Big( L^{\mp}_{21}(u_2)  L^{\pm}_{21}(u_1) b(u \mp \frac{c}{2},P)   + \ L^{\mp}_{21}(u_2)  L^{\pm}_{22}(u_1) \bar{c_0}(u \mp \frac{c}{2},P) \Big) \rho^{\pm}(u \mp \frac{c}{2}),\nn\\ \\
&&\rho^{\pm}(u \pm \frac{c}{2})  L^{\pm}_{22}(u_1)  L^{\mp}_{21}(u_2) \nn\\
&&\qquad \quad = \Big(L^{\mp}_{21}(u_2)  L^{\pm}_{21}(u_1) {c_0}(u \mp \frac{c}{2},P) + \  L^{\mp}_{21}(u_2)  L^{\pm}_{22}(u_1) \bar{b}(u \mp \frac{c}{2},P)\Big) \ \rho^{\pm}(u \mp \frac{c}{2}),\nn\\ \lb{A.29} \\
&&\rho^{\pm}(u \pm \frac{c}{2})  L^{\pm}_{22}(u_1) L^{\mp}_{22}(u_2)= L^{\mp}_{22}(u_2)  L^{\pm}_{22}(u_1) \rho^{\pm}(u \mp \frac{c}{2}),
\ena
\end{prop}
where the dynamical $R$-matrix is expressed as
\begin{center}
\bea \lb{rappendix}
&&R^\pm(u,P)=\rho^\pm(u) \left(
\begin{array}{cccc}
1 & 0 & 0 & 0 \\
0 &b(u,P) & c(u,P) & 0 \\
0 & \bar{c}(u,P) & \bar{b}(u,P) & 0 \\
0 & 0 & 0 & 1 \\
\end{array}
\right).
\ena 
\end{center}

\end{appendix}

%\LastPageEnding

\end{document}